\documentclass[11pt,oneside,reqno]{amsart}

\usepackage{amsmath,amssymb,amsthm,textcomp,mathtools}
\usepackage{amsfonts,graphicx}
\usepackage[mathscr]{eucal}
\usepackage{color}
\usepackage{csquotes}
\usepackage{enumitem}
\allowdisplaybreaks

\numberwithin{equation}{section}

\theoremstyle{definition}

\numberwithin{equation}{section}

\newtheorem{theorem}{\bf Theorem}[section]
\newtheorem{remark}{\bf Remark}[section]
\newtheorem{proposition}{Proposition}[section]
\newtheorem{lemma}{Lemma}[section]
\newtheorem{corollary}{Corollary}[section]
\newtheorem{example}{Example}[section]
\newtheorem{definition}{Definition}[section]
\newtheoremstyle
{remarkstyle}
{}
{11pt}
{}
{}
{\bfseries}
{:}
{     }
{\thmname{#1} \thmnumber{#2} }

\theoremstyle{remarkstyle}

\begin{document}
	\title{Time-changed multiparameter Poisson processes and martingales}
	\author[Pradeep Vishwakarma]{Pradeep Vishwakarma}
	\address{P. Vishwakarma, Department of Mathematics,
		Indian Statistical Institute, Kolkata, 700108, INDIA.}
	\email{vishwakarmapr.rs@gmail.com}
	\author{Kuldeep Kumar Kataria}
	\address{K. K. Kataria, Department of Mathematics,
		 Indian Institute of Technology Bhilai, Durg, 491002, INDIA.}
	 \email{kuldeepk@iitbhilai.ac.in}

	\subjclass[2020]{Primary : 60G51, 60G60; Secondary: 60G44}
	
	\keywords{multiparameter Poisson process; additive L\'evy process; multivariate subordinator, multiparameter martingale}
	\date{\today}
	
	\maketitle
	\begin{abstract}
	We consider a multiparameter extension of the homogeneous Poisson counting process, namely, the multiparameter Poisson process (MPP). We derive its various distributional properties. Also, we consider an integral of the MPP and analyze its asymptotic distribution. Thereafter, we investigate three time-changed variants of the MPP, where the time-changing components are multivariate subordinator and inverse subordinators with both dependent and independent marginals. Later, we obtain some properties of the multiparameter martingales, which are then used to derive multiparameter martingale characterizations for the MPP and one of its time-changed variants.
	\end{abstract}
\section{Introduction}
Multiparameter random processes, also known as random fields, are a natural extension of one-parameter processes. The theory of multiparameter random processes exhibits significant interactions with several other theoretical and applied fields. Khoshnevisan \cite{Khoshnevisan2002a} explores links of various types of multiparameter Markov and L\'evy processes with real and functional analysis, as well as elements of group theory and analytic number theory. The multiparameter L\'evy processes have applications in fields such as mathematical statistics, statistical mechanics, and brain imaging (see \cite{Cao1999}). For example, the additive L\'evy processes naturally arise in the analysis of multiparameter random processes, and have been a topic of study for many years (see \cite{Iafrate2024}, \cite{Khoshnevisan2002a}). The potential theory of additive L\'evy processes is used to study the Hausdorff dimension and capacity of random fractals associated with L\'evy processes (see \cite{Khoshnevisan2005}).

The Poisson point process is an example of a continuous-time and non-negative integer-valued Markov process where the time space is the set of non-negative real numbers. There have been many studies on random processes indexed by a multidimensional time space; for example, Ivanoff and Merzbach \cite{Ivanoff2000} discussed a renewal approach for general point processes, including the Poisson process as a special case. For more information on multiparameter or set indexed random processes, we refer the reader to \cite{Khoshnevisan2002a} and \cite{Merzbach1986}.
Recently, time-changed or time-fractional random processes have been a popular area of study among researchers, as they exhibit global memory properties. In the past few years, various versions of time-changed and time fractional Poisson processes have been studied (see \cite{Beghin2009}, \cite{Beghin2010}, \cite{Leskin2003}, \cite{Meerschaert2011}). For recent studies on the time fractional Poisson random field indexed by $\mathbb{R}_+^2$, we refer the reader to \cite{Aletti2018}, \cite{Beghin2020}, and \cite{Kataria2024}. Also, in \cite{Barndorff2001}, the composition of a multivariate multiparameter L\'evy process with a multivariate one parameter subordinator is introduced and studied.

Motivated by the above developments and many existing literature on time-changed multiparameter random processes (for example, see \cite{Iafrate2024}, \cite{Pedersen2004}, \cite{Pedersens2004}), we study some multiparameter point processes. Specifically, we analyze the multiparameter Poisson process and its various time-changed variants.  The characterizations of processes studied in this paper differ from some other characterizations of L\'evy processes with multidimensional index sets, which have a different definition of independent increments (see \cite{Adler1983}). 

Following the characterization of multiparameter L\'evy processes in the sense of \cite{Pedersen2004} and \cite{Pedersens2004}, we extend the theory of one parameter Poisson processes to the multiparameter case. For $d\ge1$, we consider a non-negative integer-valued random field $\{\mathscr{N}(\textbf{t}),\ \textbf{t}\in\mathbb{R}^d_+\}$ with $d$-dimensional index parameter. The variable $\mathscr{N}(\textbf{t})$ denote the total number of random points in the rectangle $[\textbf{0},\textbf{t}]=[0,t_1]\times[0,t_2]\times\dots\times[0,t_d]$ for $\textbf{t}=(t_1,t_2,\dots,t_d)\in\mathbb{R}^d_+$. It can be used to model random spatial patterns across various applied fields. For example, in telecommunication and networking, it can potentially be used to model the number of network failure across time and components, in epidemic to model the infected individual across time and space, in computer vision to model the feature points in an image and object identification on a large screen, in material science to model the impurities in a solid, and in finance and insurance to model the number of insurance claim over time and types. In particular, one of the multiparameter Poisson process $\{\mathscr{N}(\textbf{t}),\ \textbf{t}\in\mathbb{R}^d_+\}$ that we have studied is of the form $\mathscr{N}(\textbf{t})=N_1(t_1)+N_2(t_2)+\cdots+N_d(t_d)$, $ \textbf{t}=(t_1,t_2,\dots,t_d)\in\mathbb{R}^d_+$,
where $\{N_i(t_i),\ t_i\ge0\}$, $i=1,2,\dots,d$ are independent Poisson processes. As the Poisson process is L\'evy, the process $\{\mathscr{N}(\textbf{t}),\ \textbf{t}\in\mathbb{R}^d_+\}$ is an additive L\'evy process in the sense of \cite{Khoshnevisan2002a}. It can be used to model the infection in an epidemic model over time and space. Suppose $\{N_1(t),\ t\ge0\}$ counts the number of infected individuals per square kilometre and $\{N_2(t),\ t\ge0\}$ counts the number of infected individuals migrating into the population over time. Then, the process $\{N_1(t_1)+N_2(t_2),\ (t_1,t_2)\in\mathbb{R}^2_+\}$ will count the total number of infected individuals residing in the region. Note that the first indexing parameter $t_1$ denotes area in square kilometres, and the second indexing parameter $t_2$ denotes the time scale.



In Section \ref{pre}, we recall the concept of subordinators and their inverse. 

In Section \ref{sec3}, we define a multiparameter extension of the Poisson process (MPP) with index set $\mathbb{R}^d_+$, $d\ge1$. It is a multiparameter counting process with independent and stationary increments, and its one-dimensional distribution is Poisson. Its distributional properties are studied in detail. Its conditional distribution, mean, variance and autocovariance function are obtained. Moreover, we consider the integral of MPP over a rectangle in $\mathbb{R}^d_+$, and derive an asymptotic approximation of this integral process.

In Section \ref{sec4.1}, we consider the subordination of an MPP with a multivariate subordinator whose marginals are independent stable subordinators, and derive its one-dimensional distribution and governing differential equation. In Section \ref{sec4.2}, we consider a novel two parameter Poisson process time-changed by a two parameter inverse subordinator whose marginals are not necessarily independent. In Theorem \ref{tpthm}, we derive the associated governing equation that involves a two parameter non-local operator. Moreover, in Section \ref{sec4.3}, we define a time-changed MPP with an independent multiparameter inverse stable subordinator as  given in (\ref{minvsubdef}). Its various distributional properties are studied. 

In Section \ref{sec5}, we first obtain some properties of the multiparameter martingales in the sense of \cite{Khoshnevisan2002a}. Moreover, we derive a sufficient and necessary condition (see Theorem \ref{mthm2}) under which the sum of one parameter independent processes yields a multiparameter martingale. In section \ref{sec5.1}, we obtain the multiparameter martingale characterizations for the MPP and one of its fractional variants is derived.

\section{Preliminaries and notations}\label{pre}
\subsection{One parameter stable subordinator and its inverse} A one parameter non-decreasing L\'evy process, $\{S^\alpha(t),\ t\ge0\}$, $\alpha\in(0,1)$ is called an $\alpha$-stable if its Laplace transform has the following form: $$\mathbb{E}e^{-wS^\alpha(t)}=e^{-tw^\alpha},\ w>0.$$ The first passage time process $\{L^\alpha(t),\ t\ge0\}$ defined as
\begin{equation}\label{minvsubdef}
	L^\alpha(t)=\inf\{u>0:S^\alpha(u)>t\}
\end{equation}
is called an inverse $\alpha$-stable subordinator. The Laplace transform of its one-dimensional density with respect to the time variable is given as follows (see \cite{Meerschaert2011}):
\begin{equation}\label{lapminv}
	\int_{0}^{\infty}e^{-wt}\mathrm{Pr}\{{L}^\alpha(t)\in\mathrm{d}x\}=w^{\alpha-1}e^{-w^{\alpha}x},\ w>0.
\end{equation} 
For $s\ge0$ and $t\ge0$, its mean and auto covariance function are given by (see \cite{Veillette2010})
\begin{equation}\label{invsubmean}
	\mathbb{E}L^\alpha(t)=\frac{t^\alpha}{\Gamma(1+\alpha)},\ t\ge0
\end{equation}
and 
\begin{equation}\label{covinvs}
	\mathbb{C}\mathrm{ov}(L^\alpha(s),L^\alpha(t))=\frac{1}{\alpha\Gamma^2(\alpha)}\int_{0}^{\min\{s,t\}}\big((t-x)^\alpha+(s-x)^\alpha\big)x^{\alpha-1}\,\mathrm{d}x-\frac{(st)^\alpha}{\Gamma^2(1+\alpha)},
\end{equation}
respectively. Moreover, the one-dimensional density,  $l_\alpha(x,t)=\mathrm{Pr}\{L^\alpha(t)\in\mathrm{d}x\}/\mathrm{d}x$, $x\ge0$  of the inverse $\alpha$-stable subordinator solves the following fractional differential equation (see \cite{Beghin2020}):
\begin{equation}\label{invsequ}
	\partial_t^\alpha l_\alpha(x,t)=-\partial_xl_\alpha(x,t),
\end{equation}
with $l_\alpha(0,t)=t^{-\alpha}/\Gamma(1-\alpha)$, where $\partial_t^\alpha$ is the Riemann-Liouville fractional derivative defined as (see \cite{Kilbas2006})
\begin{equation}\label{rlder}
	\partial_t^\alpha f(t)=\frac{1}{\Gamma(1-\alpha)}\frac{\partial}{\partial t}\int_{0}^{t}\frac{f(s)}{(t-s)^{\alpha}}\,\mathrm{d}s,\ \alpha\in(0,1).
\end{equation}

\subsection{Multiparameter inverse subordinator}
Suppose $\alpha_i\in(0,1)$ for each $i=1,2,\dots,d$. Let $\{L_i^{\alpha_i}(t_i),\ t_i\ge0\}$, $i=1,2,\dots,d$ be independent one parameter inverse stable subordinators. Let us consider a multivariate $d$-parameter process $\{\boldsymbol{\mathscr{L}}_{\boldsymbol{\alpha}}(\textbf{t}),\ \textbf{t}\in\mathbb{R}^d_+\}$, $\boldsymbol{\alpha}=(\alpha_1,\alpha_2,\dots,\alpha_d)$ defined as follows:
\begin{equation}\label{missdef}
	\boldsymbol{\mathscr{L}}_{\boldsymbol{\alpha}}(\textbf{t})=(L_1^{\alpha_1}(t_1),L_2^{\alpha_2}(t_2),\dots,L_d^{\alpha_d}(t_d)),\ \textbf{t}=(t_1,t_2,\dots,t_d)\in\mathbb{R}^d_+.
\end{equation}
Throughout this paper, we call $\{\boldsymbol{\mathscr{L}}_{\boldsymbol{\alpha}}(\textbf{t}),\ \textbf{t}\in\mathbb{R}^d_+\}$ the multiparameter inverse stable subordinator.

\section{Multiparameter Poisson process} \label{sec3}
Here, we define a multiparameter extension of the Poisson process, namely, the multiparameter Poisson process (MPP). First, we define a multiparameter counting process as follows:	

\begin{definition}(\textbf{Multiparameter counting process})
	A non-negative integer-valued random field  $\{\mathscr{N}(\textbf{t}),\ \textbf{t}\in\mathbb{R}^d_+\}$ will be called the multiparameter counting process if almost surely, we have\\
	\noindent (i) $\mathscr{N}(\textbf{0})=0$;\\
	\noindent (ii) $\mathscr{N}(\textbf{t})\ge0$ for all $\textbf{t}\in\mathbb{R}^d_+$;\\
	\noindent (iii) $\mathscr{N}(\textbf{s})\leq\mathscr{N}(\textbf{t})$ for any $\textbf{s}\preceq\textbf{t}$.
\end{definition}
\begin{example}
	Let $\{N_i(t_i),\ t_i\ge0\}$, $i=1,2,\dots,d$ be one parameter counting processes. The multiparameter random process $\mathscr{N}=\{\mathscr{N}(\textbf{t}),\ \textbf{t}\in\mathbb{R}^d_+\}$ defined as
	\begin{equation}\label{amcp}
		\mathscr{N}(\textbf{t})\coloneqq N_1(t_1)+N_2(t_2)+\cdots+N_d(t_d),\ \textbf{t}=(t_1,t_2,\dots,t_d)\in\mathbb{R}^d_+
	\end{equation}
	is a $d$-parameter counting process.
\end{example}

\begin{definition}(\textbf{Multiparameter Poisson process})
	Let $\boldsymbol{\Lambda}=(\lambda_1,\lambda_2,\dots,\lambda_d)\in\mathbb{R}^d_{+}$ be such that $\boldsymbol{\Lambda}\succ\textbf{0}$, and let $\{\mathscr{N}(\textbf{t}),\ \textbf{t}\in\mathbb{R}^d_+\}$  be a  multiparameter counting process. We call it the MPP with transition parameter $\boldsymbol{\Lambda}$ if\\
	\noindent (i) it has independent increments, that is, for $\textbf{t}^{(i)}\in\mathbb{R}^d_+$, $i=0,1,\dots,m$ such that $\textbf{0}=\textbf{t}^{(0)}\preceq\textbf{t}^{(1)}\preceq\dots\preceq\textbf{t}^{(m)}$, the increments $\mathscr{N}(\textbf{t}^{(1)})-\mathscr{N}(\textbf{t}^{(0)})$, $\mathscr{N}(\textbf{t}^{(2)})-\mathscr{N}(\textbf{t}^{(1)})$, $\dots$, $\mathscr{N}(\textbf{t}^{(m)})-\mathscr{N}(\textbf{t}^{(m-1)})$ are independent of each other,\\
	\noindent (ii) it has stationary increments, that is, for any $\textbf{s}$ and $\textbf{t}$ both in $\mathbb{R}^d_+$ such that $\textbf{s}\preceq\textbf{t}$, we have $\mathscr{N}(\textbf{t})-\mathscr{N}(\textbf{s})\overset{d}{=}\mathscr{N}(\textbf{t}-\textbf{s})$, where $\overset{d}{=}$ denotes the equality in distribution,\\
	\noindent (iii) for all $\textbf{t}\in\mathbb{R}^d_+$, $\mathscr{N}(\textbf{t})$ has Poisson distribution with mean $\boldsymbol{\Lambda}\cdot\textbf{t}=\lambda_1t_1+\lambda_2t_2+\dots+\lambda_dt_d$, given as
	\begin{equation}\label{def1}
		\mathrm{Pr}\{\mathscr{N}(\textbf{t})=n\}=\frac{(\boldsymbol{\Lambda}\cdot\textbf{t})^n}{n!}e^{-\boldsymbol{\Lambda}\cdot\textbf{t}},\ n\ge0.
	\end{equation}
\end{definition}

Note that the MPP is a multiparameter L\'evy process in the sense of \cite{Khoshnevisan2002a}. Indeed, it is an integer-valued multiparameter subordinator whose Laplace transform is given by
$$	\mathbb{E}e^{-\eta\mathscr{N}(\textbf{t})}=\exp(\boldsymbol{\Lambda}\cdot\textbf{t}(e^{-\eta}-1)),\ \eta>0,
$$
whose  $d$-dimensional Laplace exponent is $\mathscr{B}(\eta)=(\lambda_1(1-e^{-\eta}), \lambda_2(1-e^{-\eta}),\dots,\lambda_d(1-e^{-\eta}))$. The probability generating function (pgf) of MPP is given by
\begin{equation}\label{mpppgf}
	\mathbb{E}u^{\mathscr{N}(\textbf{t})}=\exp(\boldsymbol{\Lambda}\cdot\textbf{t}(u-1)),\ |u|\leq1,\,\textbf{t}\in\mathbb{R}^d_+.
\end{equation}
Its mean and variance are  $\mathbb{E}\mathscr{N}(\textbf{t})=\boldsymbol{\Lambda}\cdot\textbf{t}$ and $\mathbb{V}\mathrm{ar}\mathscr{N}(\textbf{t})=\boldsymbol{\Lambda}\cdot\textbf{t}$, respectively.

\begin{remark}\label{summpp}
	(i)
	Let $\{N_i(t_i),\ t_i\ge0\}$, $i=1,2,\dots,d$ be independent one parameter counting processes such that $\{\sum_{i=1}^{d}N_i(t_i),$ $ (t_1,t_2,\dots,t_d)\in\mathbb{R}^d_+\}$ is a $d$-parameter Poisson process. Then,  $\{\sum_{j=1}^{d'}N_j(t_j),\ (t_1,t_2,\dots,t_{d'})\in\mathbb{R}^{d'}_+\}$ is a $d'$-parameter Poisson process for all $1\leq d'\leq d$.\vskip3pt
	
	\noindent (ii) From \cite{Iafrate2024} and \cite{Barndorff2001}, it follows that for any MPP $\{\mathscr{N}(\textbf{t}),\ \textbf{t}\in\mathbb{R}^d_+\}$ with transition parameter $\boldsymbol{\Lambda}=(\lambda_1,\lambda_2,\dots,\lambda_d)\succ\textbf{0}$ there exist $d$-many independent Poisson processes $\{N_1(t_1),\ t_1\ge0\}$, $\{N_2(t_2),\ t_2\ge0\}$, $\dots$, $\{N_d(t_d),\ t_d\ge0\}$ with transition rates $\lambda_1$, $\lambda_2$, $\dots$, $\lambda_d$, respectively, such that
	$\mathscr{N}(\textbf{t})\overset{d}{=}\sum_{i=1}^{d}N_i(t_i)$ for all $\textbf{t}=(t_1,t_2,\dots,t_d)\in\mathbb{R}^d_+
	$. As the Poisson process is a L\'evy process, any MPP of the form (\ref{amcp}) is an additive L\'evy process in the sense of \cite{Khoshnevisan2002a}. 
	It is well known that the auto covariance of a Poisson process $\{N(t),\ t\ge0\}$ with transition rate $\lambda>0$ is $\mathbb{C}\mathrm{ov}(N(t),N(s))=\lambda\min\{s,t\}$ for all $s\ge0$ and $t\ge0$. Thus, for all $\textbf{s}$ and $\textbf{t}$ in $\mathbb{R}^d_+$, the auto covariance function of MPP is given by $\mathbb{C}\mathrm{ov}(\mathscr{N}(\textbf{s}),\mathscr{N}(\textbf{t}))=\sum_{i=1}^{d}\lambda_i\min\{s_i,t_i\}$.
\end{remark} 

Next, we show that the multiparameter counting process of the form (\ref{amcp}) is an MPP if and only if it is point wise equal to the sum of independent one parameter Poisson processes. Its proof is immediate from the above remark. We omit the details.

\begin{proposition}\label{1thm}
	Let $\{N_i(t_i),\ t_i\ge0\}$, $i=1,2,\dots,d$ be one parameter independent counting processes. Then, the multiparameter counting process $\{\mathscr{N}(\textbf{t}),\ \textbf{t}\in\mathbb{R}^d_+\}$ defined as
	\begin{equation}\label{appdef}
		\mathscr{N}(\textbf{t})=N_1(t_1)+N_2(t_2)+\dots+N_d(t_d),\ \textbf{t}=(t_1,t_2,\dots,t_d)\in\mathbb{R}^d_+
	\end{equation} 
	is a $d$-parameter Poisson process with transition parameter $\boldsymbol{\Lambda}=(\lambda_1,\lambda_2,\dots,\lambda_d)\succ\textbf{0}$ if and only if each  $\{N_i(t_i),\ t_i\ge0\}$ is a Poisson process with transition rate $\lambda_i>0$.
\end{proposition}
%

\begin{remark}
	Suppose $\mathscr{N}(t)=N_1(t)+N_2(t)+\dots+N_d(t)$, $t\ge0$, where $\{N_i(t),\ t\ge0\}$, $i=1,2,\dots,d$ are independent one parameter Poisson processes with transition rates $\lambda_i>0$. Then,  $\{\mathscr{N}(t),\ t\ge0\}$ is a one parameter Poisson process with transition rate $\lambda_1+\lambda_2+\dots+\lambda_d$. It is known as the superposition of independent Poisson processes.
	Also, the weighted sum of $\{N_i(t),\ t\ge0\}$'s, that is, $\{\sum_{i=1}^{d}iN_i(t),\ t\ge0\}$ is a generalized counting process which can perform jumps of sizes $1,2,\dots,d$ with transition rates $\lambda_1$, $\lambda_2$, $\dots$, $\lambda_d$, respectively (see \cite{Dhillon2024}). Similarly, it can be shown that the weighted sum of independent one parameter Poisson processes $\{N_i(t_i),\ t_i\ge0\}$ is a generalized multiparameter counting process.
\end{remark}

In the next result, we provide some result on the conditional distribution of MPP. Its proof follows using the properties of conditional probability, and the independent and stationary increments of MPP. We omit the details.
\begin{proposition}
	(i)
	Let $\{\mathscr{N}(\textbf{t}),\ \textbf{t}\in\mathbb{R}^d_+\}$ be an MPP with transition parameter $\boldsymbol{\Lambda}\succ\textbf{0}$. For $\textbf{s}\prec\textbf{t}$ and $m\ge1$, conditional $\{\mathscr{N}(\textbf{t})=m\}$, the one-dimensional distribution of MPP is binomial given as follows:
	\begin{equation}\label{mppcmean}
		\mathrm{Pr}\{\mathscr{N}(\textbf{s})=n|\mathscr{N}(\textbf{t})=m\}=\binom{m}{n}\left(\frac{\boldsymbol{\Lambda}\cdot\textbf{s}}{\boldsymbol{\Lambda}\cdot\textbf{t}}\right)^n\left(1-\frac{\boldsymbol{\Lambda}\cdot\textbf{s}}{\boldsymbol{\Lambda}\cdot\textbf{t}}\right)^{m-n},\ 0\leq n\leq m.
	\end{equation}
	Moreover,
	\begin{equation}\label{mppcvar}
		\mathbb{E}\{\mathscr{N}(\textbf{s})|\mathscr{N}(\textbf{t})=m\}=\frac{m\boldsymbol{\Lambda}\cdot\textbf{s}}{\boldsymbol{\Lambda}\cdot\textbf{t}}\ \ \text{and}\ \ \mathbb{V}\mathrm{ar}\{\mathscr{N}(\textbf{s})|\mathscr{N}(\textbf{t})=m\}=\frac{m\boldsymbol{\Lambda}\cdot\textbf{s}}{\boldsymbol{\Lambda}\cdot\textbf{t}}\left(1-\frac{\boldsymbol{\Lambda}\cdot\textbf{s}}{\boldsymbol{\Lambda}\cdot\textbf{t}}\right).
	\end{equation}
	\noindent (ii) For $\textbf{0}\prec\textbf{r}\preceq\textbf{s}\preceq\textbf{t}$, we have
	\begin{equation}\label{mppbcmean}
		\mathbb{E}\{\mathscr{N}(\textbf{r})\mathscr{N}(\textbf{s})|\mathscr{N}(\textbf{t})=m\}=\frac{m\boldsymbol{\Lambda}\cdot\textbf{r}}{\boldsymbol{\Lambda}\cdot\textbf{t}}+m(m-1)\frac{(\boldsymbol{\Lambda}\cdot\textbf{r})(\boldsymbol{\Lambda}\cdot\textbf{s})}{(\boldsymbol{\Lambda}\cdot\textbf{t})^2},\ m\ge1.
	\end{equation}
\end{proposition}
\begin{remark}
	Note that the expressions obtained in (\ref{mppcmean}) and (\ref{mppcvar}) depend on the transition parameter $\boldsymbol{\Lambda}$, but in the case of one parameter Poisson process, that is, for $d=1$, these are independent of the choice of transition rate. Moreover, if we take $\lambda_1=\lambda_2=\dots=\lambda_d$, then these expressions are independent of $\boldsymbol{\Lambda}$.
\end{remark}
\begin{remark}
	On using (\ref{mppbcmean}), we get
	\begin{equation*}
		\mathbb{E}\{\mathscr{N}(\textbf{r})\mathscr{N}(\textbf{s})\}=\mathbb{E}(\mathbb{E}\{\mathscr{N}(\textbf{r})\mathscr{N}(\textbf{s})|\mathscr{N}(\textbf{t})\})=(\boldsymbol{\Lambda}\cdot\textbf{r})(\boldsymbol{\Lambda}\cdot\textbf{s})+\boldsymbol{\Lambda}\cdot\textbf{r},\ \textbf{0}\prec\textbf{r}\preceq\textbf{s}.
	\end{equation*}	
	Moreover, on substituting $\textbf{r}=\textbf{s}$ in (\ref{mppbcmean}), we get
	\begin{equation*}
		\mathbb{E}\{\mathscr{N}(\textbf{s})^2|\mathscr{N}(\textbf{t})=m\}=\frac{m\boldsymbol{\Lambda}\cdot\textbf{s}}{\boldsymbol{\Lambda}\cdot\textbf{t}}+m(m-1)\bigg(\frac{\boldsymbol{\Lambda}\cdot\textbf{s}}{\boldsymbol{\Lambda}\cdot\textbf{t}}\bigg)^2,\ \textbf{s}\prec\textbf{t}.
	\end{equation*}
	Thus,
	\begin{equation*}
		\mathbb{V}\mathrm{ar}\{\mathscr{N}(\textbf{s})|\mathscr{N}(\textbf{t})=m\}=\frac{m\boldsymbol{\Lambda}\cdot\textbf{s}}{\boldsymbol{\Lambda}\cdot\textbf{t}}-m\bigg(\frac{\boldsymbol{\Lambda}\cdot\textbf{s}}{\boldsymbol{\Lambda}\cdot\textbf{t}}\bigg)^2,\ \textbf{s}\prec\textbf{t}.
	\end{equation*}
\end{remark}

The proof of the following result follows from Proposition \ref{1thm}. Hence, we omit it.
\begin{proposition}
	Let $\{\mathscr{N}(\textbf{t}),\ \textbf{t}\in\mathbb{R}^d_+\}$ be an MPP defined as $\mathscr{N}(\textbf{t})=N_1(t_1)+N_2(t_2)+\cdots+ N_d(t_d)$, $\textbf{t}=(t_1,t_2,\dots,t_d)\in\mathbb{R}^d_+$  where $\{N_1(t),\ t\ge0\}, \{N_2(t),\ t\ge0\},\dots, \{N_d(t),\ t\ge0\}$ are independent Poisson processes with rate parameters $\lambda_1,\lambda_2,\dots,\lambda_d$, respectively. The conditional distribution of $N_i(t_i)$ given $\{\mathscr{N}(\textbf{t})=m\}$, $m\ge1$ is binomial with success probability $\lambda_it_i/(\boldsymbol{\Lambda}\cdot\textbf{t})$, that is,
	\begin{equation*}
		\mathrm{Pr}\{N_i(t_i)=n|\mathscr{N}(\textbf{t})=m\}=\binom{m}{n}\bigg(\frac{\lambda_it_i}{\boldsymbol{\Lambda}\cdot\textbf{t}}\bigg)^n\bigg(\frac{\sum_{j\ne i}\lambda_jt_j}{\boldsymbol{\Lambda}\cdot\textbf{t}}\bigg)^{m-n},\  0\leq n\leq m.
	\end{equation*}
	Its conditional expectation is given by
	$$\mathbb{E}\{N_i(t_i)|\mathscr{N}(\textbf{t})=m\}=\frac{m\lambda_it_i}{\boldsymbol{\Lambda}\cdot\textbf{t}},\  \textbf{t}\in\mathbb{R}^d_+.$$ 
\end{proposition}

\subsection{Thinning of an MPP} In one parameter case, the thinning of a Poisson process is done using independent Bernoulli trials. That is, whenever there is a Poisson arrival, it is retained with probability $p\in(0,1)$ and discarded with probability $1-p$. It is well known that the resultant processes from thinning are independent Poisson processes with reduced rate parameters. Furthermore, in this case, the associated semigroup operator admits an explicit representation in terms of the thinned Poisson process, yielding a Mehler-type formula (see \cite{Last2017}, Section 20.1).

Here, we discuss the thinning of an MPP. Let $\{Y_j\}_{j\ge1}$ be a sequence of independent and identically distributed Bernoulli random variables with success probability $p\in(0,1)$. Let $\{\mathscr{N}(\textbf{t}),\ \textbf{t}\in\mathbb{R}^d_+\}$ be an MPP with rate parameter $\boldsymbol{\Lambda}$, and it is independent of $\{Y_j\}_{j\ge1}$. We consider a multiparameter process defined as
\begin{equation*}
	\mathscr{N}_1(\textbf{t})\coloneqq\sum_{j=1}^{\mathscr{N}(\textbf{t})}Y_j,\ \textbf{t}\in\mathbb{R}^d_+.
\end{equation*}
Note that $\{\mathscr{N}_1(\textbf{t}),\ \textbf{t}\in\mathbb{R}^d_+\}$ is a multiparameter version of the compound Poisson process with Bernoulli compounding variables. We observe that $\{\mathscr{N}_1(\textbf{t}),\ \textbf{t}\in\mathbb{R}^d_+\}$ is a multiparameter L\'evy process. Indeed, by an elementary computation, it can be shown that its characteristic function coincides with that of an MPP with transition parameter $p\boldsymbol{\Lambda}$. Further, we define a different multiparameter process,
\begin{equation*}
	\mathscr{N}_2(\textbf{t})\coloneqq \mathscr{N}(\textbf{t})-\mathscr{N}_1(\textbf{t}),\ \textbf{t}\in\mathbb{R}^d_+.
\end{equation*}
Using similar arguments, $\{\mathscr{N}_2(\textbf{t}),\ \textbf{t}\in\mathbb{R}^d_+\}$ is also an MPP with transition parameter $(1-p)\boldsymbol{\Lambda}$. 

\subsection{Integrals of MPP}
In \cite{Orsingher2013}, some integrals of homogeneous and fractional  Poisson processes are introduced and studied. An integral of a fractional Poisson random field over a rectangle is studied in \cite{Kataria2024}. We now consider the fractional integrals of an MPP over a rectangle in $\mathbb{R}^d_+$. The motivation for studying the integrals of random processes comes from the appearance of integrated processes in many applied mathematical literature, for example, see \cite{Puri1966}, \cite{Vishwakarma2024a}, \cite{Vishwakarma2024b}, and references therein.

For $\textbf{t}=(t_1,t_2,\dots,t_d)\in\mathbb{R}^d_+$ and $\boldsymbol{\rho}=(\rho_1,\rho_2,\dots,\rho_d)\succ\textbf{0}$, we define the Riemann-Liouville fractional integral of the MPP as follows:
\begin{equation}\label{mppint}
	\mathscr{X}^{\boldsymbol{\rho}}(\textbf{t})=\frac{1}{\Gamma(\rho_1)\dots\Gamma(\rho_d)}\int_{0}^{t_1}(t_1-s_1)^{\rho_1-1}\dots\int_{0}^{t_d}(t_d-s_d)^{\rho_d-1}\mathscr{N}(\textbf{s})\,\mathrm{d}s_1\dots\mathrm{d}s_d.
\end{equation} 
Its mean is 
\begin{equation*}
	\mathbb{E}\mathscr{X}^{\boldsymbol{\rho}}(\textbf{t})=\sum_{j=1}^{d}\frac{\lambda_jt_j^{\rho_j+1}}{\Gamma(\rho_j+2)}\prod_{i\ne j}\frac{t_i^{\rho_i}}{\Gamma(\rho_i+1)},\ \textbf{t}\in\mathbb{R}^d_+.
\end{equation*}
Also, its conditional mean is given by
\begin{align*}
	\mathbb{E}&\{\mathscr{X}^{\boldsymbol{\rho}}(\textbf{t})|\mathscr{N}(\textbf{t})=m\}\\
	&=\frac{1}{\prod_{i=1}^{d}\Gamma(\rho_i)}\int_{0}^{t_1}(t_1-s_1)^{\rho_1-1}\dots\int_{0}^{t_d}(t_d-s_d)^{\rho_d-1}\mathbb{E}\{\mathscr{N}(\textbf{s})|\mathscr{N}(\textbf{t})=m\}\,\prod_{i=1}^{d}\mathrm{d}s_i\\
	&=\frac{1}{\prod_{i=1}^{d}\Gamma(\rho_i)}\int_{0}^{t_1}(t_1-s_1)^{\rho_1-1}\dots\int_{0}^{t_d}(t_d-s_d)^{\rho_d-1}\frac{m\boldsymbol{\Lambda}\cdot\textbf{s}}{\boldsymbol{\Lambda}\cdot\textbf{t}}\,\prod_{i=1}^{d}\mathrm{d}s_i\\
	&=\frac{m}{\boldsymbol{\Lambda}\cdot\textbf{t}}\sum_{j=1}^{d}\frac{\lambda_jt_j^{\rho_j+1}}{\Gamma(\rho_j+2)}\prod_{i\ne j}\frac{t_i^{\rho_i}}{\Gamma(\rho_i+1)},\ \textbf{t}\in\mathbb{R}^d_+,
\end{align*}
where we have used $(\ref{mppcvar})$ to obtain the second equality.
\begin{remark}
	For $d=1$, $\alpha\in(0,1)$ and $\rho>0$, (\ref{mppint}) reduces to the fractional integral of  a homogeneous Poisson process (see \cite{Orsingher2013}), $$X^\rho(t)=\Gamma^{-1}(\rho)\int_{0}^{t}(t-s)^{\rho-1}N(s)\,\mathrm{d}s,$$
	whose mean and variance are  $\mathbb{E}X^\rho(t)={\lambda t^{\rho+1}}/{\Gamma(\rho+2)}$ and
	\begin{equation}\label{varppint}
		\mathbb{V}\mathrm{ar}X^\rho(t)=\frac{\lambda t^{2\rho+1}}{(2\rho+1)\Gamma^2(\rho+1)},\ t\ge0,
	\end{equation}
	respectively.
\end{remark}

If the MPP in (\ref{mppint}) is of the form (\ref{appdef}) then from Proposition \ref{1thm}, it follows that there exist independent Poisson processes $\{N_1(t_1),\ t_1\ge0\}$, $\{N_2(t_2), t_2\ge0\}$, $\dots$, $\{N_d(t_d),\ t_d\ge0\}$ with transition rates $\lambda_1$, $\lambda_2$, $\dots$, $\lambda_d$, respectively, such that 
\begin{align*}
	\mathscr{X}^{\boldsymbol{\rho}}(\textbf{t})&=\frac{1}{\Gamma(\rho_1)\dots\Gamma(\rho_d)}\int_{0}^{t_1}(t_1-s_1)^{\rho_1-1}\dots\int_{0}^{t_d}(t_d-s_d)^{\rho_d-1}\sum_{j=1}^{d}N_j(s_j)\,\mathrm{d}s_1\dots\mathrm{d}s_d\\
	&=\sum_{j=1}^{d}\prod_{i\ne j}\frac{t_i^{\rho_i}}{\Gamma(\rho_i+1)}\frac{1}{\Gamma(\rho_j)}\int_{0}^{t_j}(t_j-s_j)^{\rho_j-1}N_j(s_j)\,\mathrm{d}s_j,\ \textbf{t}\in\mathbb{R}^d_+.
\end{align*}	
On using (\ref{varppint}), we get the variance of (\ref{mppint}) as follows:
\begin{equation}\label{intvar}
	\mathbb{V}\mathrm{ar}\mathscr{X}^{\boldsymbol{\rho}}(\textbf{t})=\sum_{j=1}^{d}\frac{\lambda_j t_j^{2\rho_j+1}}{(2\rho_j+1)\Gamma^2(\rho_j+1)}\bigg(\prod_{i\ne j}\frac{t_i^{\rho_i}}{\Gamma(\rho_i+1)}\bigg)^2,\ \textbf{t}\in\mathbb{R}^d_+.
\end{equation}
In particular, for $d=1$, (\ref{intvar}) reduces to (\ref{varppint}).

The next result provides a random sum representation of the Riemann integral of MPP. 
\begin{proposition}
	Let $\{\mathscr{N}(\textbf{t}),\ \textbf{t}\in\mathbb{R}^d_+\}$ be an MPP with transition parameter $\boldsymbol{\Lambda}$ as defined in (\ref{appdef}). Then, its integral has the following representation:
	\begin{align*}
		\mathscr{X}(\textbf{t})&=\int_{0}^{t_1}\int_{0}^{t_2}\dots\int_{0}^{t_d}\mathscr{N}(\textbf{s})\,\mathrm{d}s_1\mathrm{d}s_2\dots\mathrm{d}s_d\\
		&\overset{d}{=}\sum_{j=1}^{d}\sum_{k=1}^{N_j(t_j)}Y_{jk}\prod_{i\ne j}t_i,\  \textbf{t}=(t_1,t_2,\dots,t_d)\in\mathbb{R}^d_+,
	\end{align*}
	where $Y_{j1}$, $Y_{j2}$, $Y_{j3},\dots$ are independent and identically distributed uniform $[0,t_j]$ random variables for every $j=1,2,\dots,d$. In particular, the integral of MPP is a weighted sum of independent compound Poisson processes.
\end{proposition}
\begin{proof}
	In view of Proposition (\ref{1thm}), $\{N_i(t_i),\ t_i\ge0\}$, $i=1,2,\dots,d$ are independent one parameter Poisson processes. So,
	\begin{equation*}
		\mathscr{X}(\textbf{t})=\sum_{j=1}^{d}\int_{0}^{t_1}\int_{0}^{t_2}\dots\int_{0}^{t_d}N_j(s_j)\,\mathrm{d}s_1\mathrm{d}s_2\dots\mathrm{d}s_d=\sum_{j=1}^{d}\prod_{i\ne j}t_i\int_{0}^{t_j}N_j(s_j)\,\mathrm{d}s_j.
	\end{equation*}
	For each $j\in\{1,2,\dots,d\}$ and $t_j>0$, from Theorem 4.1 of \cite{Orsingher2013}, it follows that there exist independent random variables $Y_{j1}$, $Y_{j2}$, $\dots$ with common uniform distribution over $[0,t_j]$ such that $\int_{0}^{t_j}N_j(s)\,\mathrm{d}s\overset{d}{=}\sum_{k=1}^{N_j(t_j)}Y_{jk}$. This completes the proof.
	
\end{proof}
\begin{remark}
	From Eq. (4.10) of \cite{Orsingher2013}, the characteristic function of the integral $X(t)=\int_{0}^{t}N(s)\,\mathrm{d}s$ of a Poisson process $\{N(t),\ t\ge0\}$ is given by \begin{equation}\label{ppintch}
		\mathbb{E}e^{i\eta X(t)}=\exp\bigg(\frac{i\lambda\eta t^2}{2}-\frac{\lambda\eta^2t^3}{6}+o(t^3)\bigg),\ \eta\in\mathbb{R},
	\end{equation} 
	where $o(t^3)/t^3\to0$ as $t\to0$. Thus, the characteristic function of $\mathscr{X}(\textbf{t})$ is given by
	\begin{align*}
		\mathbb{E}e^{i\eta\mathscr{X}(\textbf{t})}&=\mathbb{E}\exp\bigg(i\eta\sum_{j=1}^{d}\prod_{l\ne j}t_l\int_{0}^{t_j}N_j(s_j)\,\mathrm{d}s_j\bigg)\\
		&=\prod_{j=1}^{d}\mathbb{E}\exp\bigg(i\eta\prod_{l\ne j}t_l\int_{0}^{t_j}N_j(s_j)\,\mathrm{d}s_j\bigg)\\
		&=\prod_{j=1}^{d}\exp\bigg(\frac{i\lambda_j\eta\prod_{l\ne j}t_l t_j^2}{2}-\frac{\lambda_j(\eta\prod_{l\ne j}t_l)^2t_j^3}{6}+o(t_j^3)\bigg)\\
		&=\exp\bigg(i\eta\sum_{j=1}^{d}\frac{\lambda_j\prod_{l\ne j}t_l t_j^2}{2}-\frac{\eta^2}{2}\sum_{j=1}^{d}\frac{\lambda_j\prod_{l\ne j}t_l^2t_j^3}{3}+\sum_{j=1}^{d}o(t_j^3)\bigg),
	\end{align*}
	where we have used (\ref{ppintch}) to get the penultimate step. Thus, for a $\textbf{t}\in\mathbb{R}^d_+$ with sufficiently small components, the integrated MPP has a Gaussian distribution with mean $\sum_{j=1}^{d}{\lambda_j\prod_{l\ne j}t_l t_j^2}/{2}$ and variance $\sum_{j=1}^{d}{\lambda_j\prod_{l\ne j}t_l^2t_j^3}/{3}$. Moreover, the mean and variance of $\mathscr{X}(\textbf{t})$ coincide with the parameters of an approximating multiparameter Gaussian process.
\end{remark}

\section{Time-changed MPPs}\label{sec4}
Here, we investigate some time-changed variants of the MPP, where the time-changing components are a subordinator and its inverse.


\subsection{Time-changed by stable subordinators}\label{sec4.1} Let $\{\mathscr{N}(\textbf{t}),\ \textbf{t}\in\mathbb{R}^d_+\}$ be an MPP with transition parameter $\boldsymbol{\Lambda}=(\lambda_1,\lambda_2$, $\dots,\lambda_d)\succ\textbf{0}$. For $i=1,2,\dots,d$ and $\alpha_i\in(0,1)$, let $\{\boldsymbol{\mathcal{S}}_{\boldsymbol{\alpha}}(t),\ t\ge0\}$ be a one parameter multivariate subordinator defined as $\boldsymbol{\mathcal{S}}_{\boldsymbol{\alpha}}(t)\coloneqq(S_1^{\alpha_1}(t),S_2^{\alpha_2}(t)$, $\dots,S_d^{\alpha_d}(t))$. Its marginals are independent stable subordinators which are independent of $\{\mathscr{N}(\textbf{t}),\ \textbf{t}\in\mathbb{R}^d_+\}$. We define a time-changed multiparameter process as follows:
\begin{equation}\label{tcss}
	\mathscr{N}^{\boldsymbol{\alpha}}(t)\coloneqq\mathscr{N}(\boldsymbol{\mathcal{S}}_{\boldsymbol{\alpha}}(t)),\ \boldsymbol{\alpha}=(\alpha_1,\alpha_2,\dots,\alpha_d),\ t\ge0.
\end{equation}
It is a one parameter L\'evy process whose Laplace transform is given by $\mathbb{E}e^{-\eta\mathscr{N}^{\boldsymbol{\alpha}}(t)}=\exp(-t\sum_{i=1}^{d}$ $\lambda_i^{\alpha_i}(1-e^{-\eta})^{\alpha_i})$, $\eta>0$. From Remark \ref{summpp}, it follows that  $\mathscr{N}^{\boldsymbol{\alpha}}(t)\overset{d}{=}N_1(S_1^{\alpha_1}(t))+$ $N_2(S_2^{\alpha_2}(t))+\dots+N_1(S_d^{\alpha_d}(t))$, where $\{N_i(t),\ t\ge0\}$'s are independent Poisson processes. For $d=1$, the process defined by (\ref{tcss}) reduces to the space fractional Poisson process studied in \cite{Orsingher2012}. 
\begin{remark}
	In view of Eq. (2.21) of \cite{Orsingher2012}, the process $\{\mathscr{N}^{\boldsymbol{\alpha}}(t),\ t\ge0\}$ is equal in distribution to a finite sum of independent space fractional Poisson processes. For each $i=1,2,\dots,d$, the distribution $p^{\alpha_i}_i(n,t)=\mathrm{Pr}\{N_i(S_i^{\alpha_i}(t))=n\}$, $n\ge0$ is given by
	\begin{equation}\label{sfppdist}
		p^{\alpha_i}_i(n,t)=\frac{(-1)^n}{n!}\sum_{r=0}^{\infty}\frac{(-\lambda_i^{\alpha_i}t)^r\Gamma(\alpha_ir+1)}{r!\Gamma(\alpha_ir+1-n)},
	\end{equation}
	which solves the following system of differential equations:
	\begin{equation}\label{sfppdiff}
		\frac{\mathrm{d}}{\mathrm{d}t}p^{\alpha_i}_i(n,t)=-\lambda_i^{\alpha_i}(I-B)^{\alpha_i}p^{\alpha_i}_i(n,t)=-\lambda_i^{\alpha_i}\sum_{r=0}^{n}\frac{(-1)^r\Gamma(\alpha_i+1)}{r!\Gamma(\alpha_i+1-r)}p^{\alpha_i}_i(n-r,t),\ n\ge0,
	\end{equation}
	with initial conditions $p^{\alpha_i}_i(0,0)=1$ and $p^{\alpha_i}_i(n,0)=0$ for all $n\ge1$, and also $p^{\alpha_i}_i(n,t)=0$ whenever $n<0$. Here, $B$ denotes the backward shift operator, that is, $B(p^{\alpha_i}_i(n,t))=p^{\alpha_i}_i(n-1,t)$ and
	\begin{equation*}
		(I-B)^\alpha=\sum_{r=0}^{\infty}(-1)^r\binom{\alpha}{r}B^r.
	\end{equation*}
\end{remark}
The system of equations (\ref{sfppdiff}) corresponds to the governing equation of the state probabilities of a Poisson process time-changed by an independent subordinator, which results in a L\'evy process. As established in \cite{Orsingher2012}, it can be considered as a generalization of the Kolmogorov forward and backward equations of the Poisson process, obtained by replacing $I-B$ by a fractional difference operator $(I-B)^\alpha$, $0<\alpha\leq 1$. In particular, on letting $\alpha\rightarrow1$, these equations reduce to the Kolmogorov equations associated with the homogeneous Poisson process.

Next, we obtain the distribution of (\ref{tcss}) and its governing system of differential equations.
\begin{theorem}\label{thmsub}
	The distribution $p^{\boldsymbol{\alpha}}(n,t)=\mathrm{Pr}\{\mathscr{N}^{\boldsymbol{\alpha}}(t)=n\}$, $n\ge0$  is given by
	\begin{equation*}
		p^{\boldsymbol{\alpha}}(n,t)=\sum_{\Theta(n,d)}\prod_{i=1}^{d}\frac{(-1)^{n_i}}{n_i!}\sum_{r=0}^{\infty}\frac{(-\lambda_i^{\alpha_i}t)^r\Gamma(\alpha_ir+1)}{r!\Gamma(\alpha_ir+1-n_i)},
	\end{equation*}
	where  $\Theta(n,d)=$ $\{(n_1,n_2,\dots,n_d): 0\leq n_i\leq n,\, \sum_{i=1}^{d}n_i=n\}$. It solves the following system of differential equations:
	\begin{equation*}
		\frac{\mathrm{d}}{\mathrm{d}t}p^{\boldsymbol{\alpha}}(n,t)=-\sum_{j=1}^{d}\lambda_j^{\alpha_j}(I-B)^{\alpha_j}p^{\boldsymbol{\alpha}}(n,t),\ n\ge0,
	\end{equation*}
	with initial conditions $p^{\boldsymbol{\alpha}}(0,0)=1$ and $p^{\boldsymbol{\alpha}}(n,0)=0$ for all $n\ge1$. 
\end{theorem}
\begin{proof}
	We have $p^{\boldsymbol{\alpha}}(n,t)=\sum_{\Theta(n,d)}\prod_{i=1}^{d}p_i^{\alpha_i}(n_i,t)$, $n\ge0$, $t\ge0$. Thus, using (\ref{sfppdist}) yields the required distribution. 
	
	On taking the derivative of $p^{\boldsymbol{\alpha}}(n,t)$ with respect to $t$, we get
	\begin{align*}
		\frac{\mathrm{d}}{\mathrm{d}t}p^{\boldsymbol{\alpha}}(n,t)&=\sum_{\Theta(n,d)}\sum_{j=1}^{d}\prod_{i\ne j}p_i^{\alpha_i}(n_i,t)\frac{\mathrm{d}}{\mathrm{d}t}p_j^{\alpha_j}(n_j,t)\\
		&=-\sum_{\Theta(n,d)}\sum_{j=1}^{d}\lambda_j^{\alpha_j}\sum_{r_j=0}^{n_j}\frac{(-1)^{r_j}\Gamma(\alpha_j+1)}{r_j!\Gamma(\alpha_j+1-r_j)}\prod_{i\ne j}p_i^{\alpha_i}(n_i,t)p^{\alpha_j}_j(n_j-r_j,t)\\
		&=-\sum_{j=1}^{d}\lambda_j^{\alpha_j}\sum_{r_j=0}^{n}\frac{(-1)^{r_j}\Gamma(\alpha_j+1)}{r_j!\Gamma(\alpha_j+1-r_j)}\sum_{\Theta(n-r_j,d)}\prod_{i=1}^{d}p_i^{\alpha_i}(n_i,t)\\
		&=-\sum_{j=1}^{d}\lambda_j^{\alpha_j}\sum_{r_j=0}^{\infty}\frac{(-1)^{r_j}\Gamma(\alpha_j+1)}{r_j!\Gamma(\alpha_j+1-r_j)}p^{\boldsymbol{\alpha}}(n-r_j,t)\\
		&=-\sum_{j=1}^{d}\lambda_j^{\alpha_j}\sum_{r_j=0}^{\infty}\frac{(-1)^{r_j}\Gamma(\alpha_j+1)}{r_j!\Gamma(\alpha_j+1-r_j)}B^{r_j}p^{\boldsymbol{\alpha}}(n,t).
	\end{align*}
	This completes the proof.
\end{proof}
\begin{remark}
	The probability generating function $G^{\boldsymbol{\alpha}}(u,t)=\mathbb{E}u^{\mathscr{N}^{\boldsymbol{\alpha}}(t)}$, $|u|\leq1$ solves 
	\begin{align*}
		\frac{\mathrm{d}}{\mathrm{d}t}G^{\boldsymbol{\alpha}}(u,t)&=-\sum_{j=1}^{d}\lambda_j^{\alpha_j}\sum_{n=0}^{\infty}u^n\sum_{r_j=0}^{n}\frac{(-1)^{r_j}\Gamma(\alpha_j+1)}{r_j!\Gamma(\alpha_j+1-r_j)}p^{\boldsymbol{\alpha}}(n-r_j,t)\\
		&=-\sum_{j=1}^{d}\lambda_j^{\alpha_j}\sum_{r_j=0}^{\infty}\frac{(-u)^{r_j}\Gamma(\alpha_j+1)}{r_j!\Gamma(\alpha_j+1-r_j)}\sum_{n=r_j}^{\infty}u^{n-r_j}p^{\boldsymbol{\alpha}}(n-r_j,t)\\
		&=-\sum_{j=1}^{d}\lambda_j^{\alpha_j}(1-u)^{\alpha_j}G^{\boldsymbol{\alpha}}(u,t),
	\end{align*}
	with $G^{\boldsymbol{\alpha}}(u,0)=1$. It is given by $$G^{\boldsymbol{\alpha}}(u,t)=\exp\bigg(-\sum_{j=1}^{d}\lambda_j^{\alpha_j}(1-u)^{\alpha_j}t\bigg),\ |u|\leq1.$$ 
	Hence, $\mathbb{E}(\mathscr{N}^{\boldsymbol{\alpha}}(t))^k=\infty$ for all $k\ge1$.
\end{remark}

\subsection{MPP time-changed by two parameter inverse subordinator}\label{sec4.2} We now consider the case of two parameter Poisson process, that is, the MPP with $d=2$. We study its composition with a bivariate two parameter inverse subordinator with dependent marginals. We are considering a two parameter case only because the study of inverse subordinator with dependent marginals and more than two parameters can be quite involved.

Let $\{(S_1(t),S_2(t)),\ t\ge0\}$ be a bivariate subordinator and $\{\boldsymbol{\mathscr{L}}(t_1,t_2),\ (t_1,t_2)\in\mathbb{R}^2_+\}$ be the two parameter inverse subordinator such that $\boldsymbol{\mathscr{L}}(t_1,t_2)=(L_1(t_1),L_2(t_2))$, whose marginals $L_i(t_i)=\inf\{u>0:S_i(u)>t_i\}$, $i=1,2$ are not necessarily independent. It is observed that the density of $\boldsymbol{\mathscr{L}}(t_1,t_2)$ has two components (see \cite{Beghin2020}), the first component $l(x_1,x_2,t_1,t_2)\,\mathrm{d}x_1\mathrm{d}x_2=\mathrm{Pr}\{L_1(t_1)\in\mathrm{d}x_1,L_2(t_2)\in\mathrm{d}x_2\}$, $x_1\ne x_2$ is  absolutely continuous and the second component $\tilde{l}(x,t_1,t_2)\,\mathrm{d}x=\mathrm{Pr}\{L_1(t_1)\in\mathrm{d}x,L_2(t_2)\in\mathrm{d}x\}$ has support on the line $x_1=x_2$.

We need the following result (see \cite{Beghin2020}, Theorem 3.6):
\begin{theorem}
	The absolutely continuous part  $l(x_1,x_2,t_1,t_2)$ solves 
	\begin{equation}\label{biinvs1dist}
		\mathcal{D}_{t_1,t_2}l(x_1,x_2,t_1,t_2)=-(\partial_{x_1}+\partial_{x_2})l(x_1,x_2,t_1,t_2)
	\end{equation}
	when $x_2>x_1>0$, with boundary condition
	\begin{equation*}
		l(0,x_2,t_1,t_2)=(\mathcal{D}_{t_1,t_2}-\mathcal{D}_{t_2}^{(2)})\mathcal{D}_{t_2}^{(2)}\mathcal{H}(t_1)\mathrm{Pr}\{L_2(t_2)\ge x_2\},\ x_2>0,
	\end{equation*}
	and for $x_1>x_2>0$, it solves (\ref{biinvs1dist}) under the boundary condition
	\begin{equation*}
		l(x_1,0,t_1,t_2)=(\mathcal{D}_{t_1,t_2}-\mathcal{D}_{t_1}^{(1)})\mathcal{D}_{t_1}^{(1)}\mathcal{H}(t_2)\mathrm{Pr}\{L_1(t_1)\ge x_1\},\ x_1>0.
	\end{equation*}
	Here, $\mathcal{H}$ is the Heaviside function and $\mathcal{D}^{(i)}_t$ is an operator defined by
	\begin{equation}
		\mathcal{D}^{(i)}_tg(t)\coloneqq\int_{0}^{\infty}(g(t)-g(t-s))\nu_i(\mathrm{d}s),
	\end{equation}
	where $\nu_i$ is the L\'evy measure of $S_i(t)$, $i=1,2$.
	
	Further, $\tilde{l}(x,t_1,t_2)$, $x\ge0$ solves
	\begin{equation}\label{biinvs2dist}
		\mathcal{D}_{t_1,t_2}\tilde{l}(x,t_1,t_2)=-\partial_x\tilde{l}(x,t_1,t_2),
	\end{equation}
	with boundary condition 
	\begin{equation*}
		\tilde{l}(0,t_1,t_2)=\bar{\nu}(t_1,t_2)=\int_{t_1}^{\infty}\int_{t_2}^{\infty}\nu(\mathrm{d}y_1,\mathrm{d}y_2),
	\end{equation*}
	where the operator $\mathcal{D}_{t_1,t_2}$ is defined as 
	\begin{equation*}
		\mathcal{D}_{t_1,t_2}f(x_1,x_2)=\int_{0}^{\infty}\int_{0}^{\infty}\left(f(x_1,x_2)-f(x_1-y_1,x_2-y_2)\right)\nu(\mathrm{d}y_1,\mathrm{d}y_2),
	\end{equation*} 
	and $\nu$ is the L\'evy measure associated with the bivariate subordinator $(S_1(t),S_2(t))$.
\end{theorem}

Let $\{\mathscr{N}(t_1,t_2),\ (t_1,t_2)\in\mathbb{R}^2_+\}$ be a two parameter Poisson process with transition parameter $(\lambda_1,\lambda_2)\succ(0,0)$. Its one-dimensional distribution $p(n,t_1,t_2)=\mathrm{Pr}\{\mathscr{N}(t_1,t_2)=n\}$, $n\ge0$ solves the following forward equations:
\begin{equation}\label{bippequ}
	(\partial_{t_1}+\partial_{t_2})p(n,t_1,t_2)=-(\lambda_1+\lambda_2)(I-B)p(n,t_1,t_2),\ n\ge0.
\end{equation} 
For $i=1,2$, we have
\begin{equation}\label{parbippequ}
	\partial_{t_i}p(n,t_1,t_2)=-\lambda_i(I-B)p(n,t_1,t_2),\ n\ge0.
\end{equation}

We now consider a time-changed two parameter process $\{\mathscr{N}^*(t_1,t_2),\ (t_1,t_2)\in\mathbb{R}^2_+\}$ defined as follows:
\begin{equation}\label{bpsubmpp}
	\mathscr{N}^*(t_1,t_2)=\mathscr{N}(L_1(t_1),L_2(t_2)),\ (t_1,t_2)\in\mathbb{R}^2_+.
\end{equation}

Next result gives a system of governing differential equations for the marginal distribution of the time-changed process (\ref{bpsubmpp}). Unlike the bivariate time-changed Markov process studied in Theorem 4.1 of \cite{Beghin2020}, these equations involve the initial conditions of the density of two parameter inverse subordinator because  $\mathrm{Pr}\{\mathscr{N}(0,t_2)=n\}$ is not necessarily zero for $n\ne0$.
\begin{theorem}\label{tpthm}
	The distribution $p^*(n,t_1,t_2)=\mathrm{Pr}\{\mathscr{N}^*(t_1,t_2)=n\}$, $n\ge0$ of (\ref{bpsubmpp}) satisfies the following system of equations:
	\begin{align*}
		\mathcal{D}_{t_1,t_2}p^*(n,t_1,t_2)&=-(\lambda_1+\lambda_2)(I-B)p^*(n,t_1,t_2)\nonumber\\
		&\ \ +\frac{(-\lambda_2)^n}{n!}\partial_{\lambda_2}^n\lambda_2^{-1}(\mathcal{D}_{t_1,t_2}-\mathcal{D}_{t_2}^{(2)})\mathcal{D}_{t_2}^{(2)}\mathcal{H}(t_1)(1-\mathbb{E}e^{-\lambda_2L_2(t_2)})\nonumber\\
		&\ \ +\frac{(-\lambda_1)^n}{n!}\partial_{\lambda_1}^n\lambda_1^{-1}(\mathcal{D}_{t_1,t_2}-\mathcal{D}_{t_1}^{(1)})\mathcal{D}_{t_1}^{(1)}\mathcal{H}(t_2)(1-\mathbb{E}e^{-\lambda_1L_1(t_1)}),\ n\ge1,
	\end{align*}
	and 
	\begin{align*}
		\mathcal{D}_{t_1,t_2}p^*(0,t_1,t_2)&=-(\lambda_1+\lambda_2)p^*(0,t_1,t_2)\nonumber\\
		&\ \ +\lambda_2^{-1}(\mathcal{D}_{t_1,t_2}-\mathcal{D}_{t_2}^{(2)})\mathcal{D}_{t_2}^{(2)}\mathcal{H}(t_1)(1-\mathbb{E}e^{-\lambda_2L_2(t_2)})\nonumber\\
		&\ \ +\lambda_1^{-1}(\mathcal{D}_{t_1,t_2}-\mathcal{D}_{t_1}^{(1)})\mathcal{D}_{t_1}^{(1)}\mathcal{H}(t_2)(1-\mathbb{E}e^{-\lambda_1L_1(t_1)})+\bar{\nu}(\mathrm{d}t_1,\mathrm{d}t_2),
	\end{align*}
	with initial condition $p^*(0,0,0)=1$, where $Bp^*(n,t_1,t_2)=p^*(n-1,t_1,t_2)$ and $\bar{\nu}(t_1,t_2)=\int_{t_1}^{\infty}\int_{t_2}^{\infty}\nu(\mathrm{d}s_1,\mathrm{d}s_2)$.
\end{theorem}
\begin{proof}
	The one-dimensional distribution of $\{\mathscr{N}^*(t_1,t_2),\ (t_1,t_2)\in\mathbb{R}^2_+\}$ is given by
	\begin{align}\label{bimppdist}
		p^*(n,t_1,t_2)&=\int_{0}^{\infty}\int_{0}^{\infty}p(n,x_1,x_2)l(x_1,x_2,t_1,t_2)\,\mathrm{d}x_1\,\mathrm{d}x_2\nonumber\\
		&\ \ +\int_{0}^{\infty}p(n,x,x)\tilde{l}(x,t_1,t_2)\,\mathrm{d}x,
	\end{align}
	where $p(n,t_1,t_2)$ is the distribution of two parameter Poisson process. 
	
	On applying the operator $\mathcal{D}_{t_1,t_2}$ on both sides of (\ref{bimppdist}) and using the fact that it commutes with the integral, we get
	\begin{align*}
		\mathcal{D}_{t_1,t_2}p^*(n,t_1,t_2)&=-\int_{0}^{\infty}\int_{0}^{\infty}p(n,x_1,x_2)\partial_{x_1}l(x_1,x_2,t_1,t_2)\,\mathrm{d}x_1\,\mathrm{d}x_2\\
		&\ \ -\int_{0}^{\infty}\int_{0}^{\infty}p(n,x_1,x_2)\partial_{x_2}l(x_1,x_2,t_1,t_2)\,\mathrm{d}x_1\,\mathrm{d}x_2\\
		&\ \ -\int_{0}^{\infty}p(n,x,x)\partial_{x}\tilde{l}(x,t_1,t_2)\,\mathrm{d}x,
	\end{align*}
	where we have used (\ref{biinvs1dist}) and (\ref{biinvs2dist}). Integrating by part and using (\ref{bippequ})-(\ref{parbippequ}) along with $\mathscr{N}(0,0)=0$, we get
	\begin{align*}
			\mathcal{D}_{t_1,t_2}p^*(n,t_1,t_2)
			&=-\lambda_1(I-B)\int_{0}^{\infty}\int_{0}^{\infty}p(n,x_1,x_2)l(x_1,x_2,t_1,t_2)\,\mathrm{d}x_1\,\mathrm{d}x_2\\
			&\ \ -\lambda_2(I-B)\int_{0}^{\infty}\int_{0}^{\infty}p(n,x_1,x_2)l(x_1,x_2,t_1,t_2)\,\mathrm{d}x_1\,\mathrm{d}x_2\\
			&\ \ +\int_{0}^{\infty}p(n,0,x_2)l(0,x_2,t_1,t_2)\,\mathrm{d}x_2+\int_{0}^{\infty}p(n,x_1,0)l(x_1,0,t_1,t_2)\,\mathrm{d}x_1\\
			&\ \ -(\lambda_1+\lambda_2)(I-B)\int_{0}^{\infty}p(n,x,x)\tilde{l}(x,t_1,t_2)\,\mathrm{d}x+p(n,0,0)\bar{\nu}(\mathrm{d}t_1,\mathrm{d}t_2)\\
			&=-(\lambda_1+\lambda_2)(I-B)p^*(n,t_1,t_2)+\frac{(-\lambda_2)^n}{n!}\partial_{\lambda_2}^n\int_{0}^{\infty}e^{-\lambda_2x_2}l(0,x_2,t_1,t_2)\,\mathrm{d}x_2\\
			&\ \ +\frac{(-\lambda_1)^n}{n!}\partial_{\lambda_1}^n\int_{0}^{\infty}e^{-\lambda_1x_1}l(x_1,0,t_1,t_2)\,\mathrm{d}x_1+\delta_0(n)\bar{\nu}(\mathrm{d}t_1,\mathrm{d}t_2)\\
			&=-(\lambda_1+\lambda_2)(I-B)p^*(n,t_1,t_2)\\
			&\ \ +\frac{(-\lambda_2)^n}{n!}\partial_{\lambda_2}^n\lambda_2^{-1}(\mathcal{D}_{t_1,t_2}-\mathcal{D}_{t_2}^{(2)})\mathcal{D}_{t_2}^{(2)}\mathcal{H}(t_1)(1-\mathbb{E}e^{-\lambda_2L_2(t_2)})\\
			&\ \ +\frac{(-\lambda_1)^n}{n!}\partial_{\lambda_1}^n\lambda_1^{-1}(\mathcal{D}_{t_1,t_2}-\mathcal{D}_{t_1}^{(1)})\mathcal{D}_{t_1}^{(1)}\mathcal{H}(t_2)(1-\mathbb{E}e^{-\lambda_1L_1(t_1)})+\delta_0(n)\bar{\nu}(\mathrm{d}t_1,\mathrm{d}t_2),
	\end{align*}
	where we have used 
	\begin{equation*}
		\lambda_i^{-1}(1-\mathbb{E}e^{-\lambda_iL_i(t_i)})=\int_{0}^{\infty}e^{-\lambda_ix_i}\mathrm{Pr}\{L_i(t_i)\ge x_i\}\,\mathrm{d}x_i,\ \ i=1,2.
	\end{equation*}
	This completes the proof.
\end{proof}

\subsection{MPP time-changed by inverse stable subordinators}\label{sec4.3}
Suppose $\alpha_i\in(0,1)$ for $i=1,2,,\dots,d$. Let $\{\boldsymbol{\mathscr{L}}_{\boldsymbol{\alpha}}(\textbf{t}),\ \textbf{t}\in\mathbb{R}^d_+\}$ be a multiparameter inverse stable sobordinator as defined in (\ref{missdef}), that is, $\boldsymbol{\mathscr{L}}_{\boldsymbol{\alpha}}(\textbf{t})=(L_1^{\alpha_1}(t_1),L_2^{\alpha_2}(t_2),\dots,L_d^{\alpha}(t_d))$ with independent marginals. Let $\{\mathscr{N}(\textbf{t}),\ \textbf{t}\in\mathbb{R}^d_+\}$ be an MPP with transition parameter $\boldsymbol{\Lambda}=(\lambda_1,\lambda_2,\dots,\lambda_d)\succ\textbf{0}$ that is independent of $\{\boldsymbol{\mathscr{L}}_{\boldsymbol{\alpha}}(\textbf{t}),\ \textbf{t}\in\mathbb{R}^d_+\}$. We consider a multiparameter process $\{\mathscr{N}_{\boldsymbol{\alpha}}(\textbf{t}),\ \textbf{t}\in\mathbb{R}^d_+\}$ defined as follows: 
\begin{equation}\label{mfppdef}
	\mathscr{N}_{\boldsymbol{\alpha}}(\textbf{t})\coloneqq\mathscr{N}(\boldsymbol{\mathscr{L}}_{\boldsymbol{\alpha}}(\textbf{t})),\ \textbf{t}\in\mathbb{R}^d_+.
\end{equation}
We call it the multiparameter fractional Poisson process (MFPP) with transition parameter $\boldsymbol{\Lambda}$. For $d=1$, it reduces to a one parameter fractional Poisson process in the sense of \cite{Meerschaert2011}.

Similar to Proposition \ref{1thm}, the following result shows that a multiparameter counting process of the form (\ref{amcp}) is an MFPP if and only if it is equal to the sum of independent one parameter fractional Poisson processes.
\begin{theorem}\label{summfppthm}
	Let $\{N_1^{\alpha_1}(t_1),\ t_1\ge0\}$, $\{N_2^{\alpha_2}(t_2),\ t_2\ge0\}$, $\dots$, $\{N_d^{\alpha_d}(t_d),\ t_d\ge0\}$ be independent one parameter counting processes, where $\alpha_i\in(0,1)$ for all $i=1,2,\dots,d$. A multiparameter counting process $\{\mathscr{N}_{\boldsymbol{\alpha}}(\textbf{t}),\ \textbf{t}\in\mathbb{R}^d_+\}$ defined by $\mathscr{N}_{\boldsymbol{\alpha}}(\textbf{t})=\sum_{i=1}^{d}N_i^{\alpha_i}(t_i)$, $\textbf{t}=(t_1,t_2,\dots,t_d)\in\mathbb{R}^d_+$ is  an MFPP with transition parameter $\boldsymbol{\Lambda}=(\lambda_1,\lambda_2,\dots,\lambda_d)\succ\textbf{0}$ if and only if  $\{N_i^{\alpha_i}(t_i),\ t_i\ge0\}$ is a one parameter fractional Poisson process with transition rate $\lambda_i>0$ for each $i=1,2,\dots,d$.
\end{theorem} 
\begin{proof}
	If $\{N_i^{\alpha_i}(t_i),\ t_i\ge0\}$ is a fractional Poison process with positive transition rate $\lambda_i$ for each $i=1,2,\dots,d$, then there exist a Poisson process $\{N_i(t_i),\ t_i\ge0\}$ with transition rate $\lambda_i$, and an independent inverse stable subordinator $\{L_i^{\alpha_i}(t_i),\ t_i\ge0\}$ such that $N_i^{\alpha_i}(t_i)=N_i(L_i^{\alpha_i}(t_i))$ for all $t_i\ge0$. From Proposition \ref{1thm}, it follows that $\{\sum_{i=1}^{d}N_i(t_i),\ (t_1,t_2,\dots,t_d)\in\mathbb{R}^d_+\}$ is an MPP with transition parameter $\boldsymbol{\Lambda}=(\lambda_1,\lambda_2,\dots,\lambda_d)\succ\textbf{0}$. Thus, the process $\{\sum_{i=1}^{d}N_i(L_i^{\alpha_i}(t_i)),\ (t_1,t_2,\dots,t_d)\in\mathbb{R}^d_+\}$ is an MFPP.
	
	Conversely, if $\sum_{i=1}^{d}N_i^{\alpha_i}(t_i)$ is an MFPP with transition parameter $\boldsymbol{\Lambda}\succ\textbf{0}$ then there exist an MPP, $\{\mathscr{N}(\textbf{t}),\ \textbf{t}\in\mathbb{R}^d_+\}$, and an independent multiparameter inverse stable subordinator $\{\boldsymbol{\mathscr{L}}_{\boldsymbol{\alpha}}(\textbf{t}),\ \textbf{t}\in\mathbb{R}^d_+\}$ as defined in (\ref{missdef}) such that $\sum_{i=1}^{d}N_i^{\alpha_i}(t_i)=\mathscr{N}(\boldsymbol{\mathscr{L}}_{\boldsymbol{\alpha}}(\textbf{t}))$ for all $\textbf{t}=(t_1,t_2,\dots,t_d)\in\mathbb{R}^d_+$. Take $\textbf{t}(i)=(0,\dots,0,t_i,0,\dots,0)\in\mathbb{R}^d_+$ for $i=1,2,\dots,d$. Then, $N_i^{\alpha_i}(t_i)=\mathscr{N}(\boldsymbol{\mathscr{L}}_{\boldsymbol{\alpha}}(\textbf{t}(i)))$ almost surely, which is a composition of $N_i(t_i)=\mathscr{N}(\textbf{t}(i))$, $t_i\ge0$ and $L_i^{\alpha_i}(t_i)$. From the definition of MPP, it follows that $\{N_i(t_i),\ t_i\ge0\}$ is a Poisson process with transition rate $\lambda_i>0$. Thus, $\{N_i^{\alpha_i}(t_i),\ t_i\ge0\}$ is a fractional Poisson process for each $i=1,2,\dots,d$. This completes the proof.
\end{proof}

\begin{theorem}
	The distribution of MFPP is given by
	\begin{equation}\label{mfppdist}
		p_{\boldsymbol{\alpha}}(n,\textbf{t})=\mathrm{Pr}\{\mathscr{N}_{\boldsymbol{\alpha}}(\textbf{t})=n\}=
		\sum_{\Theta(n,d)}\prod_{i=1}^{d}(\lambda_it_i^{\alpha_i})^{n_i}E_{\alpha_i,n_i\alpha_i+1}^{n_i+1}(-\lambda_i t_i^{\alpha_i}),\ n\ge0,
	\end{equation}
	where the set $\Theta(n,d)$ is as defined in Theorem \ref{thmsub}. Here, $E_{\alpha,\beta}^\gamma$, $\alpha>0$, $\beta>0$, $\gamma>0$ is the three parameter Mittag-Leffler function defined as follows (see \cite{Kilbas2006}):
	\begin{equation*}\label{3mittag}
		E_{\alpha,\beta}^\gamma(x)=\sum_{k=0}^{\infty}\frac{(\gamma)_kx^k}{\Gamma(k\alpha+\beta)k!},\ x\in\mathbb{R},
	\end{equation*}
	where $(\gamma)_k=\gamma(\gamma+1)\dots(\gamma+k-1)$.
	
	Moreover, for each $i=1,2,\dots,d$, (\ref{mfppdist}) solves the following system of fractional differential equations:
	\begin{equation}\label{idfe}
		\mathcal{D}_{t_i}^{\alpha_i}p_{\boldsymbol{\alpha}}(n,\textbf{t})=-\lambda_i(p_{\boldsymbol{\alpha}}(n,\textbf{t})-p_{\boldsymbol{\alpha}}(n-1,\textbf{t})),\ n\ge0,
	\end{equation}
	where  $\mathcal{D}_{t_i}^{\alpha_i}$ is the Caputo fractional derivative defined by
	\begin{equation*}
		\mathcal{D}_t^\alpha f(t)=\partial_t^\alpha f(t)-\frac{t^{-\alpha}}{\Gamma(1-\alpha)}f(0^+),\ \alpha\in(0,1).
	\end{equation*}
	Here, $\partial_t^\alpha$ is the Riemann-Liouville fractional derivative as defined in (\ref{rlder}).
\end{theorem}
\begin{proof}
	Let $l_{\alpha_i}(x,t)$, $x\ge0$ be the one-dimensional density of $\{L_i^{\alpha_i}(t),\ t\ge\}$ for $i=1,2,\dots,d$. Then, in view of (\ref{mfppdef}), the distribution of MFPP is given by   
	\begin{equation}\label{mfppdistpf1}
		p_{\boldsymbol{\alpha}}(n,\textbf{t})=\int_{\mathbb{R}^d_+}p(n,\textbf{x})\prod_{i=1}^{d}l_{\alpha_i}(x_i,t_i)\in\mathrm{d}x_i\},
	\end{equation} 
	whose multivariate Laplace transform is
	\begin{align}
		\int_{\mathbb{R}^d_+}e^{-\textbf{w}\cdot\textbf{t}}p_{\boldsymbol{\alpha}}(n,\textbf{t})\,\mathrm{d}\textbf{t}
		&=\int_{\mathbb{R}^d_+}\frac{(\boldsymbol{\Lambda}\cdot\textbf{x})^n}{n!}e^{-\boldsymbol{\Lambda}\cdot\textbf{x}}\prod_{i=1}^{d}
		w_i^{\alpha_i-1}e^{-w_i^{\alpha_i}x_i}\,\mathrm{d}x_i\nonumber\\
		&=\frac{\prod_{i=1}^{d}
			w_i^{\alpha_i-1}}{n!}\sum_{\Theta(n,d)}\frac{n!}{n_1!n_2!\dots n_d!}\prod_{i=1}^{d}\lambda_i^{n_i}\int_{0}^{\infty}x_i^{n_i}e^{-(\lambda_i+w_i^{\alpha_i})}\,\mathrm{d}x_i\nonumber\\
		&=\sum_{\Theta(n,d)}\prod_{i=1}^{d}\frac{\lambda_i^{n_i}w_i^{\alpha_i-1}}{(w_i^{\alpha_i}+\lambda_i)^{n_i+1}},\ \textbf{w}=(w_1,w_2,\dots,w_d)\succ\textbf{0},\label{mdist1}
	\end{align}
	where we have used (\ref{lapminv}) and (\ref{def1}) to get the first equality. On taking the inverse Laplace transform on both sides of (\ref{mdist1}) and using Eq. (1.9.13) of \cite{Kilbas2006}, we get the required distribution.
	
	On applying the Riemann-Liouville fractional derivative operator $\partial_{t_i}^{\alpha_i}$ on both sides of (\ref{mfppdistpf1}) and using (\ref{invsequ}), we have
	\begin{align}
		\partial_{t_i}^{\alpha_i}p_{\boldsymbol{\alpha}}(n,\textbf{t})&=-\int_{\mathbb{R}^d_+}p_{\boldsymbol{\alpha}}(n,\textbf{x})\partial_{x_i}l_{\alpha_i}(x_i,t_i)\,\mathrm{d}x_i\prod_{i\ne j}l_{\alpha_j}(x_j,t_j)\,\mathrm{d}x_j\nonumber\\
		&=-\int_{\mathbb{R}^d_+}\partial_{x_i}(p_{\boldsymbol{\alpha}}(n,\textbf{x})l_{\alpha_i}(x_i,t_i))\,\mathrm{d}x_i\prod_{i\ne j}l_{\alpha_j}(x_j,t_j)\,\mathrm{d}x_j\nonumber\\ &\ \ +\int_{\mathbb{R}^d_+}\partial_{x_i}p_{\boldsymbol{\alpha}}(n,\textbf{x})l_{\alpha_i}(x_i,t_i)\,\mathrm{d}x_i\prod_{i\ne j}l_{\alpha_j}(x_j,t_j)\,\mathrm{d}x_j\nonumber\\
		&=\frac{t_i^{-\alpha_i}}{\Gamma(1-\alpha_i)}\int_{\mathbb{R}^{d-1}_+}p_{\boldsymbol{\alpha}}(n,\textbf{x}^i(0))\prod_{i\ne j}l_{\alpha_j}(x_j,t_j)\,\mathrm{d}x_j\nonumber\\ &\ \ +\int_{\mathbb{R}^d_+}\partial_{x_i}p_{\boldsymbol{\alpha}}(n,\textbf{x})l_{\alpha_i}(x_i,t_i)\,\mathrm{d}x_i\prod_{i\ne j}l_{\alpha_j}(x_j,t_j)\,\mathrm{d}x_j,\label{cpf0}
	\end{align}
	where $\textbf{x}^i(0)=(x_1,x_2,\dots,x_{i-1},0,x_{i+1},\dots,x_d)\in\mathbb{R}^d_+$. On using
	$$
	\partial_{t_i}p(n,\textbf{t})=-\lambda_i(p(n,\textbf{t})-p(n-1,\textbf{t})),\ n\ge0$$ 
	for each $i=1,2,\dots,d$
	and
	\begin{equation*}
		\int_{\mathbb{R}^{d-1}_+}p_{\boldsymbol{\alpha}}(n,\textbf{x}^i(0))\prod_{i\ne j}l_{\alpha_j}(x_j,t_j)\,\mathrm{d}x_j=p_{\boldsymbol{\alpha}}(n,\textbf{t}^i(0))
	\end{equation*}
	in (\ref{cpf0}), we get
	\begin{equation*}
		\partial_{t_i}^{\alpha_i}p_{\boldsymbol{\alpha}}(n,\textbf{t})=-\lambda_i(p_{\boldsymbol{\alpha}}(n,\textbf{t})-p_{\boldsymbol{\alpha}}(n-1,\textbf{t}))+\frac{t_i^{-\alpha_i}}{\Gamma(1-\alpha_i)}p_{\boldsymbol{\alpha}}(n,\textbf{t}^i(0)).
	\end{equation*}
	This completes the proof.
\end{proof}
\begin{remark}
	In one parameter case, that is, $d=1$, the governing system of differential equations (\ref{idfe}) reduces to the system of governing equations for the state probabilities of a fractional Poisson process introduced in \cite{Meerschaert2011}. These equations are called the Kolmogorov fractional integro-differential equations in the sense of \cite{Orsingher2018}. Consequently, the system of equations (\ref{idfe}) can be regarded as a multiparameter analogue of these equations in the sense of \cite{Iafrate2024}. 
\end{remark}

\begin{remark}
	The regularity of $p_{\boldsymbol{\alpha}}(n,t)$ can be verified as follows:
	\begin{align*}
		\int_{\mathbb{R}^d_+}e^{-\textbf{w}\cdot\textbf{t}}\sum_{n=0}^{\infty}p_{\boldsymbol{\alpha}}(n,\textbf{t})\,\mathrm{d}\textbf{t}&=\sum_{n=0}^{\infty}
		\sum_{\Theta(n,d)}\prod_{i=1}^{d}\frac{\lambda_i^{n_i}w_i^{\alpha_i-1}}{(w_i^{\alpha_i}+\lambda_i)^{n_i+1}},\ \textbf{w}=(w_1,w_2,\dots,w_d)\succ\textbf{0},\\
		&=\sum_{n_1=0}^{\infty}\sum_{n_2=0}^{\infty}\dots\sum_{n_d=0}^{\infty}\prod_{i=1}^{d}\frac{w_i^{\alpha_i-1}\lambda_i^{n_i}}{(w_i^{\alpha_i}+\lambda_i)^{n_i+1}}\\
		&=\prod_{i=1}^{d}\frac{w_i^{\alpha_i-1}}{w_i^{\alpha_i}+\lambda_i}\frac{w_i^{\alpha_i}+\lambda_i}{w_i^{\alpha_i}}=\prod_{i=1}^{d}\frac{1}{w_i},
	\end{align*} 
	whose inverse Laplace transform is $\sum_{n=0}^{\infty}p_{\boldsymbol{\alpha}}(n,\textbf{t})=1$.
\end{remark}
\begin{remark}\label{mfpprep} From (\ref{lapminv}), we have $\mathbb{E}e^{-w_iL_i^{\alpha_i}(t_i)}=E_{\alpha_i,1}(-\lambda_it_i^{\alpha_i})$, $t_i>0$, where $E_{\alpha,1}$, $\alpha>0$ is the one parameter Mittag-Leffler function. From
	Theorem 2.2 of \cite{Meerschaert2011}, it follows that $\mathscr{N}_{\boldsymbol{\alpha}}(\textbf{t})=\sum_{i=1}^{d}N_i^{\alpha_i}(t_i)$, is a finite sum of independent renewal processes with Mittag-Leffler inter-arrival times. If $\{\mathscr{N}_{\boldsymbol{\alpha}}(\textbf{t}),\ \textbf{t}\in\mathbb{R}^d_+\}$ is an MFPP with transition parameter $\boldsymbol{\Lambda}=(\lambda_1,\lambda_2,\dots,\lambda_d)\succ\textbf{0}$ then 
	$\mathscr{N}_{\boldsymbol{\alpha}}(\textbf{t})\overset{d}{=}\sum_{i=1}^{d}N_i({L}_i^{\alpha_i}(t_i))$, $\textbf{t}\in\mathbb{R}^d_+$, 
	where for each $i=1,2,\dots,d$, $\{N_i(t_i),\ t_i\ge0\}$ is a Poisson process with transition rates $\lambda_i$ and $\{L_i^{\alpha_i}(t_i),\ t_i\ge0\}$ is an inverse stable subordinator with index $\alpha_i\in(0,1)$.
	All the processes appearing here are independent of each other. Thus, we conclude that every MFPP  is equal in distribution with a finite sum of independent one parameter fractional Poisson processes. 
	
\end{remark}

\begin{remark}
	For $d=1$, the MFPP reduces to the fractional Poisson process $\{N_\alpha(t),\ t\ge0\}$, $\alpha\in(0,1)$ with transition rate $\lambda>0$ (see \cite{Beghin2009}, \cite{Meerschaert2011}). Its distribution is given by $p_\alpha(n,t)=(\lambda t^{\alpha})^nE_{\alpha,n\alpha+1}^{n+1}(-\lambda t^\alpha)$, $n\ge0$. Also, its mean and variance are 
	\begin{equation}\label{fppmean}
		\mathbb{E}N_\alpha(t)=\frac{\lambda t^\alpha}{\Gamma(\alpha+1)}\ \mathrm{and}\ 
		\mathbb{V}\mathrm{ar}N_\alpha(t)=\frac{\lambda t^\alpha}{\Gamma(\alpha+1)}+\frac{(\lambda t^{\alpha})^2}{\alpha}\left(\frac{1}{\Gamma(2\alpha)}-\frac{1}{\alpha\Gamma^2(\alpha)}\right),
	\end{equation}
	for all $t\ge0$, respectively.
\end{remark}

From Remark \ref{mfpprep}, it follows that the MFPP $\{\mathscr{N}_{\boldsymbol{\alpha}}(\textbf{t}),\ \textbf{t}\in\mathbb{R}^d_+\}$ is equal in distribution to the sum of $d$-many independent fractional Poisson processes. Thus, on using (\ref{fppmean}), we get 
\begin{equation*}\label{afppmean}
	\mathbb{E}\mathscr{N}_{\boldsymbol{\alpha}}(\textbf{t})=\sum_{i=1}^{d}\frac{\lambda_i t_i^{\alpha_i}}{\Gamma(\alpha_i+1)},\ \textbf{t}=(t_1,t_2,\dots,t_d)\in\mathbb{R}^d_+
\end{equation*}
and 
\begin{equation}\label{mfppvar}
	\mathbb{V}\mathrm{ar}\mathscr{N}_{\boldsymbol{\alpha}}(\textbf{t})=\sum_{i=1}^{d}\left(\frac{\lambda_i t_i^{\alpha_i}}{\Gamma(\alpha_i+1)}+\frac{(\lambda_i t_i^{\alpha_i})^2}{\alpha_i}\left(\frac{1}{\Gamma(2\alpha_i)}-\frac{1}{\alpha_i\Gamma^2(\alpha_i)}\right)\right),\ \textbf{t}\in\mathbb{R}^d_+.
\end{equation}

\section{Martingale characterizations}\label{sec5}
In this section, we give the martingale characterizations for  MPP and MFPP. First, we recall the definition of a multiparameter martingale in the sense of \cite{Khoshnevisan2002a}. 
Let $(\Omega,\mathscr{F},\mathbb{P})$ be a complete probability space. A collection $\{\mathscr{F}(\textbf{t}),\ \textbf{t}\in\mathbb{R}^d_+\}$  of complete sub-sigma fields of $\mathscr{F}$ is called a $d$-parameter filtration if $\mathscr{F}(\textbf{s})\subseteq\mathscr{F}(\textbf{t})$ for all $\textbf{s}\preceq \textbf{t}$.

\begin{definition}
	A $\mathbb{P}$-integrable  random process $\{M(\textbf{t}),\ \textbf{t}\in\mathbb{R}^d_+\}$ is called $d$-parameter martingale if $M(\textbf{t})$ is $\mathscr{F}(\textbf{t})$-measurable for all $\textbf{t}\in\mathbb{R}^d_+$ and $\mathbb{E}(M(\textbf{t})|\mathscr{F}(\textbf{s}))=M(\textbf{s})$ almost surely whenever $\textbf{s}\preceq\textbf{t}$.
\end{definition}
For more details on some other characterizations of multiparameter martingales, we refer the reader to \cite{Nualart1987} and \cite{Zakai1981}. 

For an identical indexing set, it is known that the sum and product of one parameter independent martingales are again one parameter martingale with respect to their natural filtrations. 

Next, we prove that the sum and the product of one parameter martingales with different indexing sets are multiparameter martingales with respect to an appropriate choice of filtration.
\begin{theorem}\label{mthm1}
	Let $\{M_i(t_i),\ t_i\ge0\}$, $i=1,2,\dots,d$ be independent one parameter martingales with respect to natural filtration $\mathscr{F}^{(i)}(t_i)=\sigma(M_i(s_i),\,s_i\leq t_i)$. Then, the collection $\{\mathscr{F}(\textbf{t}),\ \textbf{t}\in\mathbb{R}^d_+\}$ defined by $\mathscr{F}(\textbf{t})=\bigvee_{i=1}^{d}\mathscr{F}^{(i)}(t_i)$, $\textbf{t}=(t_1,t_2,\dots,t_d)$ is a multiparameter filtration and the processes $\{\sum_{i=1}^{d}M_i(t_i),\ \textbf{t}\in\mathbb{R}^d_+\}$ and $\{\prod_{i=1}^{d}M_i(t_i),\ \textbf{t}\in\mathbb{R}^d_+\}$ are multiparameter $\mathscr{F}(\textbf{t})$-martingales. 
\end{theorem}
\begin{proof}
	See \ref{app1}.
\end{proof}
\begin{theorem}\label{mthm3}
	Let $\{M_1(\textbf{t}),\ \textbf{t}\in\mathbb{R}^d_+\}$ and $\{M_2(\textbf{t}),\ \textbf{t}\in\mathbb{R}^d_+\}$ be two independent multiparameter martingales with respect to their natural filtrations $\{\mathscr{F}_1(\textbf{t}),\ \textbf{t}\in\mathbb{R}^d_+\}$ and $\{\mathscr{F}_2(\textbf{t}),\ \textbf{t}\in\mathbb{R}^d_+\}$, respectively, where $\mathscr{F}_i(\textbf{t})=\sigma(M_i(\textbf{s}),\ \textbf{s}\preceq\textbf{t})$, $i=1,2$. Then, $\{M_1(\textbf{t})+M_2(\textbf{t}),\ \textbf{t}\in\mathbb{R}^d_+\}$ and $\{M_1(\textbf{t})M_2(\textbf{t}),\ \textbf{t}\in\mathbb{R}^d_+\}$ are multiparameter martingales with respect to the filtration $\{\mathscr{F}(\textbf{t}),\ \textbf{t}\in\mathbb{R}^d_+\}$ such that $\mathscr{F}(\textbf{t})=\mathscr{F}_1(\textbf{t})\vee\mathscr{F}_2(\textbf{t})$.
\end{theorem}
\begin{proof}
	See \ref{app2}.
\end{proof}
In the following result, we show that the converse of the first part of Theorem \ref{mthm1} is true in a particular case. 
\begin{theorem}\label{mthm2}
	Let $\{M_1(t_1),\ t_1\ge0\}$, $\{M_2(t_2),\ t_2\ge0\}$, $\dots$, $\{M_d(t_d),\ t_d\ge0\}$ be independent one parameter integrable processes. For $i=1,2,\dots,d$, let $\mathscr{F}^{(i)}(t_i)=\sigma(M_i(s_i):s_i\leq t_i)$ and $\mathscr{F}(\textbf{t})=\vee_{i=1}^{d}\mathscr{F}^{(i)}(t_i)$ for all $\textbf{t}=(t_1,t_2,\dots,t_d)\in\mathbb{R}^d_+$. Then, the following statements are equivalent:\\
	\noindent (i) the process $\{\sum_{i=1}^{d}M_i(t_i),\ \textbf{t}\in\mathbb{R}^d_+\}$ is a $d$-parameter $\mathscr{F}(\textbf{t})$-martingale such that for each $i=1,2,\dots,d$, we have $\mathbb{E}M_i(s_i)=\mathbb{E}M_i(t_i)$ for all $s_i\ge0$ and $t_i\ge0$,\\
	\noindent (ii) the process $\{\sum_{i=1}^{d'}M_i(t_i),\ \textbf{t}\in\mathbb{R}^{d'}_+\}$ is a $d'$-parameter martingale with respect to $d'$-parameter filtration $\{\mathscr{F}'(\textbf{t}'),\ \textbf{t}'\in\mathbb{R}^{d'}_+\}$ for any  $1\leq d'\leq d$, where $\mathscr{F}'(\textbf{t}')=\vee_{i=1}^{d'}\mathscr{F}^{(i)}(t_i)$ for all $\textbf{t}'=(t_1,t_2,\dots,t_{d'})\in\mathbb{R}^{d'}_+$,\\
	\noindent (iii) the process $\{M_i(t_i),\ t_i\ge0\}$ is a one parameter martingale with respect to its natural filtration for each $i=1,2,\dots,d$.
\end{theorem}
\begin{proof}
	See \ref{app3}.
\end{proof}
\subsection{Martingale characterization of MPP} \label{sec5.1}
A martingale characterization of the Poisson process was given by Watanabe \cite{Watanabe1964} using its inter-arrival times. An alternate proof of this result using the stochastic calculus can be found in \cite{Shreve2004}. Here, we give a multiparameter martingale characterization for the MPP.

The following result is known as the Watanabe characterization of the Poisson process:
\begin{theorem}\label{mpp}
	Let $\{N(t),\ t\ge0\}$ be a  simple locally finite point process, and let $\{\mathscr{F}(t),\ t\ge0\}$ be a filtration such that $\mathscr{F}(t)=\sigma(N(s), s\leq t)$. Then, $\{N(t)-\lambda t,\ t\ge0\}$, $\lambda>0$ is a $\{\mathscr{F}(t),\ t\ge0\}$-martingale if and only if $\{N(t),\ t\ge0\}$ is a Poisson process with transition rate $\lambda$.
\end{theorem}

Next, we give a multiparameter martingale characterization for the MPP.
\begin{theorem}
	Let $\{N_1(t_1),\ t_1\ge0\}$, $\{N_2(t_2),\ t_2\ge0\}$, $\dots$, $\{N_d(t_d),\ t_d\ge0\}$ be simple and locally finite independent one parameter counting processes, and let $\{\mathscr{F}^{(i)}_{t_i},\ t_i\ge0\}$ be one parameter filtration such that $\mathscr{F}^{(i)}{(t_i)}=\sigma(N_i(s_i):s_i\leq t_i)$ for all $i=1,2,\dots,d$. Then, the multiparmeter counting process $\{\mathscr{N}(\textbf{t}),\ \textbf{t}\in\mathbb{R}^d_+\}$ where $\mathscr{N}(\textbf{t})=\sum_{i=1}^{d}N_i(t_i)$, $\textbf{t}=(t_1,t_2,\dots,t_d)\in\mathbb{R}^d_+$, is an MPP if and only if there exist a constant vector $\boldsymbol{\Lambda}=(\lambda_1,\lambda_2,\dots,\lambda_d)\succ\textbf{0}$  such that $\{\mathscr{N}(\textbf{t})-\boldsymbol{\Lambda}\cdot\textbf{t},\ \textbf{t}\in\mathbb{R}^d_+\}$ is a multiparameter martingale with respect to $\{\mathscr{F}(\textbf{t}),\ \textbf{t}\in\mathbb{R}^d_+\}$, where $\mathscr{F}(\textbf{t})=\bigvee_{i=1}^{d}\mathscr{F}^{(i)}{(t_i)}$ for all $\textbf{t}\in\mathbb{R}^d_+$, and $\mathbb{E}(N_i(t_i)-\lambda_it_i)=0$ for each $i=1,2,\dots,d$.
\end{theorem}
\begin{proof}
	Suppose $\{\mathscr{N}(\textbf{t}),\ \textbf{t}\in\mathbb{R}^d_+\}$ where $\mathscr{N}(\textbf{t})=\sum_{i=1}^{d}N_i(t_i)$ is a $d$-parameter Poisson process with transition parameter $\boldsymbol{\Lambda}=(\lambda_1,\lambda_2,\dots,\lambda_d)\succ\textbf{0}$. 
	If $\mathscr{F}(\textbf{t})=\bigvee_{i=1}^{d}\mathscr{F}^{(i)}({t_i})$ then $\{\mathscr{F}(\textbf{t}),\ \textbf{t}\in\mathbb{R}^d_+\}$ is a $d$-parameter filtration and $\{\mathscr{N}(\textbf{t})-\boldsymbol{\Lambda}\cdot\textbf{t},\ \textbf{t}\in\mathbb{R}^d_+\}$ is an integrable process which is adapted to $\{\mathscr{F}(\textbf{t}),\ \textbf{t}\in\mathbb{R}^d_+\}$. For $\textbf{s}\preceq\textbf{t}$ in $\mathbb{R}^d_+$, we have
	\begin{align*}
		\mathbb{E}(\mathscr{N}(\textbf{t})-\boldsymbol{\Lambda}\cdot\textbf{t}|\mathscr{F}(\textbf{s}))&=\mathbb{E}(\mathscr{N}(\textbf{t})|\mathscr{F}(\textbf{s}))-\boldsymbol{\Lambda}\cdot\textbf{t}\\
		&=\mathbb{E}(\mathscr{N}(\textbf{t})-\mathscr{N}(\textbf{s})|\mathscr{F}(\textbf{s}))+\mathscr{N}(\textbf{s})-\boldsymbol{\Lambda}\cdot\textbf{t}\\
		&=\mathbb{E}(\mathscr{N}(\textbf{t})-\mathscr{N}(\textbf{s}))+\mathscr{N}(\textbf{s})-\boldsymbol{\Lambda}\cdot\textbf{t}\\
		&=\boldsymbol{\Lambda}\cdot(\textbf{t}-\textbf{s})+\mathscr{N}(\textbf{s})-\boldsymbol{\Lambda}\cdot\textbf{t},
	\end{align*}
	where we have used the independent increments property of the MPP to get the penultimate step. Thus, $\{\mathscr{N}(\textbf{t})-\boldsymbol{\Lambda}\cdot\textbf{t},\ \textbf{t}\in\mathbb{R}^d_+\}$ is a multiparameter martingale with respect to $\{\mathscr{F}(\textbf{t}),\ \textbf{t}\in\mathbb{R}^d_+\}$. From Proposition \ref{1thm}, it follows that $\{N_i(t_i),\ t_i\ge0\}$ is a one parameter Poisson process with transition rate $\lambda_i>0$. So, $\mathbb{E}(N_i(t_i)-\lambda_it_i)=0$ for all $i=1,2,\dots,d$.
	
	Let $\{\mathscr{N}(\textbf{t})-\boldsymbol{\Lambda}\cdot\textbf{t},\ \textbf{t}\in\mathbb{R}^d_+\}$ be a multiparameter $\{\mathscr{F}(\textbf{t}),\ \textbf{t}\in\mathbb{R}^d_+\}$-martingale such that $\mathscr{N}(\textbf{t})-\boldsymbol{\Lambda}\cdot\textbf{t}=\sum_{i=1}^{d}(N_i(t_i)-\lambda_i t_i)$ and $\mathbb{E}(N_i(t_i)-\lambda_it_i)=0$ for all $i=1,2,\dots,d$. Then, from Theorem \ref{mthm2}, it follows that $\{N_1(t_1)-\lambda_1 t_1,\ t_1\ge0\}$, $\{N_2(t_2)-\lambda_2t_2,\ t_2\ge0\}$, $\dots$, $\{N_d(t_d)-\lambda_d t_d,\ t_d\ge0\}$ are independent one parameter martingales with respect to filtrations $\{\mathscr{F}^{(1)}(t_1),\ t_1\ge0\}$, $\{\mathscr{F}^{(2)}(t_2),\ t_2\ge0\}$, $\dots$, $\{\mathscr{F}^{(d)}(t_d),\ t_d\ge0\}$, respectively. This can also be directly observed by using the fact that $N_i(t_i)=N(0,\dots,0,t_i,0,\dots,0)$. So, in view of Theorem \ref{mpp},  $\{N_1(t_1),\ t_1\ge0\}$, $\{N_2(t_2),\ t_2\ge0\}$, $\dots$, $\{N_d(t_d),\ t_d\ge0\}$ are independent one parameter Poisson processes with transition rates $\lambda_1$, $\lambda_2$, $\dots$, $\lambda_d$, respectively. Thus, from Proposition \ref{1thm}, we conclude that $\{\sum_{i=1}^{d}N_i(t_i),\ \textbf{t}\in\mathbb{R}^d_+\}$ is an MPP with transition parameter $(\lambda_1,\lambda_2,\dots,\lambda_d)$. This completes the proof.
\end{proof}

\begin{remark}
	Let $\{\mathscr{F}(\textbf{t}),\ \textbf{t}\in\mathbb{R}^d_+\}$ be a multiparameter filtration, and let $\{\mathscr{N}(\textbf{t}),\ \textbf{t}\in\mathbb{R}^d_+\}$ be $\{\mathscr{F}(\textbf{t}),\ \textbf{t}\in\mathbb{R}^d_+\}$-adapted MPP with transition parameter $\boldsymbol{\Lambda}=(\lambda_1,\lambda_2,\dots,\lambda_d)$ such that $\mathscr{N}(\textbf{t})-\mathscr{N}(\textbf{s})$ is independent of $\mathscr{F}(\textbf{s})$ whenever $\textbf{s}\prec\textbf{t}$. Then, $\{\mathscr{N}(\textbf{t})-\boldsymbol{\Lambda}\cdot\textbf{t},\ \textbf{t}\in\mathbb{R}^d_+\}$ is a $d$-parameter martingale with respect to the filtration $\{\mathscr{F}(\textbf{t}),\ \textbf{t}\in\mathbb{R}^d_+\}$.
\end{remark}

The following result provides a martingale characterization for any MPP.

\begin{theorem}\label{thm4}
	The multiparamter counting process $\{\mathscr{N}(\textbf{t}),\ \textbf{t}\in\mathbb{R}^d_+\}$ is an MPP if and only if there exist a constant vector $\boldsymbol{\Lambda}=(\lambda_1,\lambda_2,\dots,\lambda_d)\succ\textbf{0}$ such that, for $c>-1$, the process $\{\exp(\mathscr{N}(\textbf{t})\ln(1+c)-c\boldsymbol{\Lambda}\cdot\textbf{t}),\ \textbf{t}\in\mathbb{R}^d_+\}$ is a multiparameter martingale with respect to the multiparameter filtration $\{\mathscr{F}(\textbf{t}),\ \textbf{t}\in\mathbb{R}^d_+\}$, where $\mathscr{F}(\textbf{t})=\sigma(\mathscr{N}(\textbf{s}),\,\textbf{s}\preceq\textbf{t})$.
\end{theorem}
\begin{proof}
	Let $\{\mathscr{N}(\textbf{t}),\ \textbf{t}\in\mathbb{R}^d_+\}$ be an MPP with transition parameter $\boldsymbol{\Lambda}$. It is enough to verify the martingale property for the process $\{\exp(\mathscr{N}(\textbf{t})\ln(1+c)-c\boldsymbol{\Lambda}\cdot\textbf{t}),\ \textbf{t}\in\mathbb{R}^d_+\}$. For $\textbf{0}\preceq \textbf{s}\preceq \textbf{t}$, by using the independent and stationary increments properties of MPP, we have
	{\small\begin{align*}
			\mathbb{E}&(\exp(\mathscr{N}(\textbf{t})\ln(1+c)-c\boldsymbol{\Lambda}\cdot\textbf{t})|\mathscr{F}(\textbf{s}))\\
			&=\mathbb{E}(\exp((\mathscr{N}(\textbf{t})-\mathscr{N}(\textbf{s}))\ln(1+c)-c\boldsymbol{\Lambda}\cdot(\textbf{t}-\textbf{s}))\exp(\mathscr{N}(\textbf{s})\ln(1+c)-c\boldsymbol{\Lambda}\cdot\textbf{s})|\mathscr{F}(\textbf{s}))\\
			&=\exp(\mathscr{N}(\textbf{s})\ln(1+c)-c\boldsymbol{\Lambda}\cdot\textbf{s})\mathbb{E}(\exp((\mathscr{N}(\textbf{t})-\mathscr{N}(\textbf{s}))\ln(1+c)-c\boldsymbol{\Lambda}\cdot(\textbf{t}-\textbf{s}))|\mathscr{F}(\textbf{s}))\\
			&=\exp(\mathscr{N}(\textbf{s})\ln(1+c)-c\boldsymbol{\Lambda}\cdot\textbf{s})\mathbb{E}\exp((\mathscr{N}(\textbf{t})-\mathscr{N}(\textbf{s}))\ln(1+c)-c\boldsymbol{\Lambda}\cdot(\textbf{t}-\textbf{s}))\\
			&=\exp(\mathscr{N}(\textbf{s})\ln(1+c)-c\boldsymbol{\Lambda}\cdot\textbf{s})\mathbb{E}(\exp(\mathscr{N}(\textbf{t}-\textbf{s})\ln(1+c)-c\boldsymbol{\Lambda}\cdot(\textbf{t}-\textbf{s})))\\
			&=\exp(\mathscr{N}(\textbf{s})\ln(1+c)-c\boldsymbol{\Lambda}\cdot\textbf{s})e^{-c\boldsymbol{\Lambda}\cdot(\textbf{t}-\textbf{s})}\exp\left(\boldsymbol{\Lambda}(\textbf{t}-\textbf{s})(e^{\ln(1+c)}-1)\right)\\
			&=\exp(\mathscr{N}(\textbf{s})\ln(1+c)-c\boldsymbol{\Lambda}\cdot\textbf{s}),
	\end{align*}}
	where the second equality follows from Theorem 34.3 of \cite{Billingsley1995} and to get the penultimate step, we have used the moment generating function of MPP, that is, $\mathbb{E}e^{x\mathscr{N}(\textbf{t})}=\exp(\boldsymbol{\Lambda}\cdot\textbf{t}(e^x-1))$, $x\in\mathbb{R}$. 
	
	For the converse part, suppose $\{\exp(\mathscr{N}(\textbf{t})\ln(1+c)-c\boldsymbol{\Lambda}\cdot\textbf{t}),\ \textbf{t}\in\mathbb{R}^d_+\}$ is a multiparameter martingale with respect to its natural filtration. Then, for any $\textbf{t}\in\mathbb{R}^d_+$ and $c=u-1$, $0<u\leq1$, we have 
	\begin{equation*}
		\mathbb{E}\exp(\mathscr{N}(\textbf{t})\ln(u)-(u-1)\boldsymbol{\Lambda}\cdot\textbf{t})=\mathbb{E}\exp(\mathscr{N}(\textbf{0})\ln(u)-(u-1)\boldsymbol{\Lambda}\cdot\textbf{0})=1.
	\end{equation*}
	Hence, $\mathbb{E}u^{\mathscr{N}(\textbf{t})}=\exp(\boldsymbol{\Lambda\cdot\textbf{t}}(u-1))$. Thus, from (\ref{mpppgf}), it follows that the distribution of $\mathscr{N}(\textbf{t})$ coincides with the distribution of an MPP. Also, for $\textbf{s}\preceq\textbf{t}$, we have 
	\begin{equation*}
		\mathbb{E}\exp\big(\mathscr{N}(\textbf{t})\ln(u)-\boldsymbol{\Lambda}\cdot\textbf{t}(u-1)|\mathscr{F}(\textbf{s})\big)=\exp\big(\mathscr{N}(\textbf{s})\ln(u)-\boldsymbol{\Lambda}\cdot\textbf{s}(u-1)\big).
	\end{equation*}
	On using Theorem 34.3 of \cite{Billingsley1995}, we have
	\begin{equation*}
		\mathbb{E}\exp\big((\mathscr{N}(\textbf{t})-\mathscr{N}(\textbf{s}))\ln(u)\big|\mathscr{F}(\textbf{s})\big)=\exp\big(\boldsymbol{\Lambda}\cdot(\textbf{t}-\textbf{s})(u-1)\big).
	\end{equation*}
	So,
	\begin{equation}\label{statpf}
		\mathbb{E}\exp\big((\mathscr{N}(\textbf{t})-\mathscr{N}(\textbf{s}))\ln(u)\big)=\exp(\boldsymbol{\Lambda}\cdot(\textbf{t}-\textbf{s})(u-1)).
	\end{equation}
	Thus, using (\ref{mpppgf}), we conclude that the process $\{\mathscr{N}(\textbf{t}),\ \textbf{t}\in\mathbb{R}^d_+\}$ has stationary increments.
	
	Let $\textbf{0}\preceq\textbf{t}^{(0)}\prec\textbf{t}^{(1)}\prec\textbf{t}^{(2)}\prec\dots\prec\textbf{t}^{(k)}$ in $\mathbb{R}^d_+$ and $0<u_l\leq1$ for all $l=1,2,\dots,k$. Then, the joint pgf of increments $\mathscr{N}(\textbf{t}^{(r+1)})-\mathscr{N}(\textbf{t}^{(r)})$, $r=0,1,\dots,k-1$ is given by
	\begin{align*}
		\mathbb{E}\prod_{r=0}^{k-1}u_{r+1}^{\mathscr{N}(\textbf{t}^{(r+1)})-\mathscr{N}(\textbf{t}^{(r)})}
		&=\mathbb{E}\exp\bigg(\sum_{r=0}^{k-1}\big(\mathscr{N}(\textbf{t}^{(r+1)})-\mathscr{N}(\textbf{t}^{(r)})\big)\ln(u_{r+1})\bigg)\\
		&=\mathbb{E}\bigg(\mathbb{E}\bigg(\exp\bigg(\sum_{r=0}^{k-1}\big(\mathscr{N}(\textbf{t}^{(r+1)})-\mathscr{N}(\textbf{t}^{(r)})\big)\ln(u_{r+1})\bigg)\Big|\mathscr{F}{(\textbf{t}^{(k-1)})}\bigg)\bigg)\\
		&=\mathbb{E}\exp\bigg(\sum_{r=0}^{k-2}\big(\mathscr{N}(\textbf{t}^{(r+1)})-\mathscr{N}(\textbf{t}^{(r)})\big)\ln(u_{r+1})\bigg)\\
		&\ \ \cdot\exp(\boldsymbol{\Lambda}\cdot(\textbf{t}^{(k)}-\textbf{t}^{(k-1)})(u_{k}-1)),\ (\text{using (\ref{statpf})})\\
		&=\mathbb{E}\bigg(\mathbb{E}\bigg(\exp\bigg(\sum_{r=0}^{k-2}\big(\mathscr{N}(\textbf{t}^{(r+1)})-\mathscr{N}(\textbf{t}^{(r)})\big)\ln(u_{r+1})\bigg)\Big|\mathscr{F}{(\textbf{t}^{(k-2)})}\bigg)\bigg)\\
		&\ \ \cdot\exp(\boldsymbol{\Lambda}\cdot(\textbf{t}^{(k)}-\textbf{t}^{(k-1)})(u_{k}-1))\\
		&\ \ \vdots\\
		&=\prod_{r=0}^{k-1}\exp(\boldsymbol{\Lambda}\cdot(\textbf{t}^{(r+1)}-\textbf{t}^{(r)})(u_{r+1}-1))=\prod_{r=0}^{k-1}\mathbb{E}u_{r+1}^{\mathscr{N}(\textbf{t}^{(r+1)})-\mathscr{N}(\textbf{t}^{(r)})},
	\end{align*}
	where the second last equality is established inductively. This proves the independent increments property of the process $\{\mathscr{N}(\textbf{t}),\ \textbf{t}\in\mathbb{R}^d_+\}$. Thus, it is an MPP with transition parameter $\boldsymbol{\Lambda}$. This completes the proof.
\end{proof}

The following result on the geometric Poisson process will be used. Two different proofs of this result using the techniques from stochastic calculus are given in \cite{Dhillon2024} and \cite{Shreve2004}. 
\begin{lemma}\label{lemma1}
	Let $\{N(t),\ t\ge0\}$ be a one parameter simple counting process with $N(0)=0$ almost surely, and let $\lambda>0$ be a positive constant and $\mathscr{F}(t)=\sigma(N(s),\ s\leq t)$. Then, for a fixed constant $c>-1$, the process $\{\exp(N(t)\ln(1+c)-\lambda ct),\ t\ge0\}$ is a martingale with respect to $\{\mathscr{F}(t),\ t\ge0\}$ whenever the compensated process $\{N(t)-\lambda t,\ t\ge0\}$ is a martingale with respect to $\{\mathscr{F}(t),\ t\ge0\}$.
\end{lemma}

The following result is a restatement of Theorem \ref{thm4} for a particular type of MPP.
\begin{corollary}
	Let $\{N_i(t_i),\ t_i\ge0\}$, $i=1,2,\dots,d$ be  independent simple counting processes, and let $\mathscr{F}^{(i)}(t_i)=\sigma(N_i(s_i),\, s_i\leq t_i)$ and $\mathscr{F}(\textbf{t})=\bigvee_{i=1}^{d}\mathscr{F}^{(i)}(t_i)$ for all $\textbf{t}=(t_1,t_2,\dots,t_d)\in\mathbb{R}^d_+$. For $\boldsymbol{\Lambda}=(\lambda_1,\lambda_2,\dots,\lambda_d)\succ\textbf{0}$, suppose $Y(\textbf{t})=\sum_{i=1}^{d}N_i(t_i)\ln(u)-\boldsymbol{\Lambda}\cdot\textbf{t}(u-1)$, $0<u\leq1$. Then, $\{\exp(Y(\textbf{t})),\ \textbf{t}\in\mathbb{R}^d_+\}$ is a multiparameter martingale with respect to $\{\mathscr{F}(\textbf{t}),\ \textbf{t}\in\mathbb{R}^d_+\}$ if and only if $\{\sum_{i=1}^{d}N_i(t_i),\ \textbf{t}\in\mathbb{R}^d_+\}$ is an MPP with transition parameter $\boldsymbol{\Lambda}$.
\end{corollary}
\begin{proof}
	If $\{\sum_{i=1}^{d}N_i(t_i),\ (t_1,t_2,\dots,t_d)\in\mathbb{R}^d_+\}$ is an MPP with transition parameter $\boldsymbol{\Lambda}$ then Proposition \ref{1thm} implies that $\{N_i(t),\ t\ge0\}$ is a Poisson process with transition rate $\lambda_i>0$ for each $i=1,2,\dots,d$. On taking $c=u-1>-1$ in Lemma \ref{lemma1}, it follows that $\{\exp(N_i(t)\ln(u)-\lambda_it(u-1)),\ t\ge0\}$'s are independent  martingales with respect to $\{\mathscr{F}^{(i)}{(t_i)},\ t\ge0\}$. Thus, in view of Theorem \ref{mthm1}, we can conclude that $\{\exp (Y(\textbf{t}))=\prod_{i=1}^{d}\exp(N_i(t_i)\ln(u)-\lambda_it_i(u-1)),\ \textbf{t}\in\mathbb{R}^d_+\}$ is a multiparameter martingale with respect to $\{\mathscr{F}(\textbf{t}),\ \textbf{t}\in\mathbb{R}^d_+\}$.
	
	The proof of the converse part follows similar lines to that of Theorem \ref{thm4}. Hence, it is omitted.
\end{proof}

\subsection{Martingale characterization of MFPP}
Aletti \textit{et al.} \cite{Aletti2018} obtained a martingale characterization for the fractional Poisson process. Here, we give a multiparameter martingale characterization of the MFPP. First, we consider the following definition of right continuity for multivariable functions (see \cite{Khoshnevisan2002a}, p. 236).
\begin{definition}
	A function $f:\mathbb{R}^d_+\to\mathbb{R}$ is right continuous if for any $\textbf{t}\in\mathbb{R}^d_+$,  
	$$\lim_{\substack{\textbf{s}\to\textbf{t}\\\textbf{t}\preceq\textbf{s}}}f(\textbf{s})=f(\textbf{t}).$$
\end{definition}
\begin{lemma}\label{rcts}
	Let $f_i:\mathbb{R}_+\to\mathbb{R}$ be  a real valued functions for each $i=1,2,\dots,d$. If $f_i$'s are right continuous on at least one point then $f:\mathbb{R}^d_+\to\mathbb{R}$ such that $f(\textbf{t})=\sum_{i=1}^{d}f_i(t_i)$, $\textbf{t}=(t_1,t_2,\dots,t_d)\in\mathbb{R}^d_+$, is a right continuous function if and only if $f_i$'s are right continuous on their respective domain.
\end{lemma}
\begin{proof}
	See \ref{app4}.
\end{proof}
\begin{theorem}\label{mfppmt}
	For $i=1,2,\dots,d$ and $\alpha_i\in(0,1)$, let $\{N_i^{\alpha_i}(t_i),\ t_i\ge0\}$ be $d$-many independent one parameter counting processes. The multiparameter  process $\{\mathscr{N}_{\boldsymbol{\alpha}}(\textbf{t})=\sum_{i=1}^{d}N_i^{\alpha_i}(t_i),\ \textbf{t}\in\mathbb{R}^d_+\}$, $\boldsymbol{\alpha}=(\alpha_1,\alpha_2,\dots,\alpha_d)$ is an MFPP if and only if there exist a constant vector $\boldsymbol{\Lambda}=(\lambda_1,\lambda_2,\dots,\lambda_d)\succ\textbf{0}$ and $d$-many independent stable subordinators $S_i^{\alpha_i}(t)$ with independent inverse processes $L_i^{\alpha_i}(t_i)=\inf\{u>0:S_i^{\alpha_i}(u)>t_i\}$ such that $\{\mathscr{N}_{\boldsymbol{\alpha}}(\textbf{t})-\boldsymbol{\Lambda}\cdot\boldsymbol{\mathscr{L}}_{\boldsymbol{\alpha}}(\textbf{t}),\ \textbf{t}\in\mathbb{R}^d_+\}$ is a right continuous multiparameter martingale with respect to $\{\mathscr{F}(\textbf{t}),\ \textbf{t}\in\mathbb{R}^d_+\}$, where $\mathscr{F}(\textbf{t})=\bigvee_{i=1}^{d}(\sigma(N_i^{\alpha_i}(s_i),\,s_i\leq t_i)\vee\sigma(L_i^{\alpha_i}(s_i),\,s_i\ge0))$. Also, $\mathbb{E}(N_i^{\alpha_i}(t_i)-\lambda_iL_i^{\alpha_i}(t_i))=0$ and for $c>0$, the collection 
	\begin{equation*}
		\{N_i^{\alpha_i}(\tau_i)-\lambda_iL_i^{\alpha_i}(\tau_i),\,\text{$\tau_i$ is a stopping time such that $L_i^{\alpha_i}(\tau_i)\leq c$}\}
	\end{equation*} 
	is uniformly integrable for each $i=1,2,\dots,d$. 
\end{theorem}
\begin{proof}
	Let $\{\mathscr{N}_{\boldsymbol{\alpha}}(\textbf{t}),\ \textbf{t}\in\mathbb{R}^d_+\}$ be an MFPP with transition parameter $\boldsymbol{\Lambda}=(\lambda_1,\lambda_2,\dots,\lambda_d)\succ\textbf{0}$. Then, $\{N_1^{\alpha_1}(t_1),\ t_1\ge0\}$, $\{N_2^{\alpha_2}(t_2),\ t_2\ge0\}$, $\dots$, $\{N_d^{\alpha_d}(t_d),\ t_d\ge0\}$ are independent fractional Poisson processes with transition rates $\lambda_1$, $\lambda_2$, $\dots$, $\lambda_d$, respectively. In view of Theorem 2 of \cite{Aletti2018}, there exist $d$-many independent stable subordinators $\{S_1^{\alpha_1}(t_1),\ t_1\ge0\}$, $\{S_2^{\alpha_2}(t_2),\ t_2\ge0\}$, $\dots$, $\{S_d^{\alpha_d}(t_d),\ t_d\ge0\}$ with corresponding inverse processes $\{L_1^{\alpha_1}(t_1),\ t_1\ge0\}$, $\{L_2^{\alpha_2}(t_2),\ t_2\ge0\}$, $\dots$, $\{L_d^{\alpha_d}(t_d),\ t_d\ge0\}$, respectively, such that $\{N_i^{\alpha_i}(t_i)-\lambda_iL_i^{\alpha_i}(t_i),\ t_i\ge0\}$ is a  right continuous martingale with respect to induced filtration $\sigma(N_i^{\alpha_i}(s_i):s_i\leq t_i)\vee\sigma(L_i^{\alpha_i}(s_i):s_i\leq t_i)$ and for $c>0$, the collection
	\begin{equation*}
		\{N_i^{\alpha_i}(\tau_i)-\lambda_iL_i^{\alpha_i}(\tau_i),\ \text{$\tau_i$ is a stopping time such that}\ L_i^{\alpha_i}(\tau_i)\leq c\}
	\end{equation*}
	is uniformly integrable	 for each $i=1,2,\dots,d$. From Theorem \ref{mthm1} and Lemma \ref{rcts}, it follows that the process $\{\sum_{i=1}^{d}(N_i^{\alpha_i}(t_i)-\lambda_iL_i^{\alpha_i}(t_i)),\ \textbf{t}\in\mathbb{R}^d_+\}$ is a right continuous multiparameter martingale with respect to $\{\mathscr{F}(\textbf{t}),\ \textbf{t}\in\mathbb{R}^d_+\}$ where $\mathscr{F}(\textbf{t})=\bigvee_{i=1}^{d}(\sigma(N_i^{\alpha_i}(s_i),\,s_i\leq t_i)\vee\sigma(L_i^{\alpha_i}(s_i),\,s_i\ge0))$ for all $\textbf{t}=(t_1,t_2,\dots,t_d)\in\mathbb{R}^d_+$. Further, using (\ref{invsubmean}) and (\ref{fppmean}), we get $\mathbb{E}(N_i^{\alpha_i}(t_i)-\lambda_iL_i^{\alpha_i}(t_i))=0$ for all $i=1,2,\dots,d$.

	Let $\boldsymbol{\mathscr{L}}_{\boldsymbol{\alpha}}(\textbf{t})=(L_1^{\alpha_1}(t_1), L_2^{\alpha_2}(t_2),\dots,L_d^{\alpha_d}(t_d))$ be a multiparameter process, where $L_i^{\alpha_i}(t_i)$'s are independent inverse stable subordinator such that $\{\sum_{i=1}^{d}N_i^{\alpha_i}(t_i)-\boldsymbol{\Lambda}\cdot\boldsymbol{\mathscr{L}}_{\boldsymbol{\alpha}}(\textbf{t}),\ \textbf{t}\in\mathbb{R}^d_+\}$ is a right continuous multiparameter martingale with respect to $\{\mathscr{F}(\textbf{t}),\ \textbf{t}\in\mathbb{R}^d_+\}$ and $\mathbb{E}(N_i^{\alpha_i}(t_i)-\lambda_iL_i^{\alpha_i}(t_i))=0$ for each $i=1,2,\dots d$. From Theorem \ref{mthm2} and Lemma \ref{rcts}, it follows that $\{N_i^{\alpha_i}(t_i)-\lambda_iL_i^{\alpha_i}(t_i),\ t_i\ge0\}$'s are independent right continuous martingales with respect to $\{\sigma(N_i^{\alpha_i}(s_i):s_i\leq t_i)\vee\sigma(L_i^{\alpha_i}(s_i):s_i\leq t_i),\ t_i\ge0\}$. Hence, $\{N_1^{\alpha_1}(t_1),\ t_1\ge0\}$, $\{N_2^{\alpha_2}(t_2),\ t_2\ge0\}$, $\dots$, $\{N_d^{\alpha_d}(t_d),\ t_d\ge0\}$ are independent fractional Poisson processes with transition rates $\lambda_1$, $\lambda_2$, $\dots$, $\lambda_d$, respectively (see \cite{Aletti2018}, Theorem 2).  Now, the proof follows from Theorem \ref{summfppthm}.
\end{proof}

\appendix
\section{Proofs of some results}
\subsection{Proof of Theorem \ref{mthm1}}\label{app1}
It is sufficient to prove the result for $d=2$. If  $\textbf{s}=(s_1,s_2)\preceq\textbf{t}=(t_1,t_2)$ in $\mathbb{R}^2_+$ then $\mathscr{F}^{(i)}(s_i)\subseteq\mathscr{F}^{(i)}(t_i)$ for $i=1,2$. Hence, $\mathscr{F}(\textbf{s})\subseteq\mathscr{F}(\textbf{t})$. Note that $M_i(t_i)$'s are $\mathscr{F}(\textbf{t})$-adapted, so is their sum and product. 

We now define the following two classes of subsets of $\mathscr{F}^{(1)}({t_1})\vee\mathscr{F}^{(2)}(t_2)$:
\begin{multline*}
	\mathcal{C}_1=\{A\cap B\in\mathscr{F}^{(1)}(s_1)\vee\mathscr{F}^{(2)}(s_2):\, A\in\mathscr{F}^{(1)}(s_1),\,B\in\mathscr{F}^{(2)}(s_2)\\ \ \text{and}\ \mathbb{E}((M_1(t_1)+M_2(t_2))I_{A\cap B})=\mathbb{E}((M_1(s_1)+M_2(s_2))I_{A\cap B})\}
\end{multline*} 
and
\begin{equation*}
	\mathcal{C}_2=\{C\in\mathscr{F}^{(1)}(s_1)\vee\mathscr{F}^{(2)}(s_2):\mathbb{E}((M_1(t_1)+M_2(t_2))I_{C})=\mathbb{E}((M_1(s_1)+M_2(s_2))I_{C})\},
\end{equation*}	
where $I_A$ denotes the indicator function on set $A$.

Note that these two classes are non-empty because for any $A\in\mathscr{F}^{(1)}(s_1)$ and $B\in\mathscr{F}^{(2)}{(s_2)}$, we have
\begin{align*}
	\mathbb{E}((M_1(t_1)+M_2(t_2))I_{A\cap B})&=\mathbb{E}(M_1(t_1)I_{A\cap B})+\mathbb{E}(M_2(t_2)I_{A\cap B})\\
	&=\mathbb{E}(M_1(t_1)I_{A})\mathbb{E}(I_{B})+\mathbb{E}(M_2(t_2)I_{B})\mathbb{E}(I_{A})\\
	&=\mathbb{E}(M_1(s_1)I_{A})\mathbb{E}(I_{B})+\mathbb{E}(M_2(s_2)I_{B})\mathbb{E}(I_{A})\\
	&=\mathbb{E}(M_1(s_1)I_{A\cap B})+\mathbb{E}(M_2(s_2)I_{A\cap B})\\
	&=\mathbb{E}((M_1(s_1)+M_2(s_2))I_{A\cap B}).
\end{align*}
Further, it can be established that $\mathcal{C}_1$ is a $\pi$-system and $\mathcal{C}_2$ is a $\lambda$-system such that $\mathcal{C}_1\subseteq\mathcal{C}_2$. So, using Dynkin's $\pi$-$\lambda$ theorem (see \cite{Billingsley1995}, Theorem 3.2), sigma algebra generated by $\mathcal{C}_1$, that is, $\sigma(\mathcal{C}_1)=\mathscr{F}^{(1)}(t_1)\vee\mathscr{F}^{(2)}(t_2)$ is contained in $\mathcal{C}_2$. Thus, $\mathcal{C}_2=\mathscr{F}^{(1)}(t_1)\vee\mathscr{F}^{(2)}(t_2)$. Therefore, $\{M_1(t_1)+M_2(t_2),\ (t_1,t_2)\in\mathbb{R}^2_+\}$ is a  $\mathscr{F}^{(1)}({t_1})\vee\mathscr{F}^{(2)}(t_2)$-martingale.

If $M_1(t_1)$ and $M_2(t_2)$ are independent then for any $A\in\mathscr{F}^{(1)}(s_1)$ and $B\in\mathscr{F}^{(2)}{(s_2)}$, we have
\begin{align*}
	\mathbb{E}(M_1(t_1)M_2(t_2)I_{A\cap B})&=\mathbb{E}(M_1(t_1)I_{A})\mathbb{E}(M_2(t_2)I_{B})\\
	&=\mathbb{E}(M_1(s_1)I_{A})\mathbb{E}(M_2(s_2)I_{B})=\mathbb{E}(M_1(s_1)M_2(s_2)I_{A\cap B}).
\end{align*}
Define two classes of subsets of $\mathscr{F}^{(1)}({t_1})\vee\mathscr{F}^{(2)}(t_2)$ as follows:
\begin{multline*}
	\mathscr{C}_1=\{A\cap B\in\mathscr{F}^{(1)}(s_1)\vee\mathscr{F}^{(2)}(s_2):\, A\in\mathscr{F}^{(1)}(s_1),\,B\in\mathscr{F}^{(2)}(s_2)\\ \ \text{and}\ \mathbb{E}(M_1(t_1)M_2(t_2)I_{A\cap B})=\mathbb{E}(M_1(s_1)M_2(s_2)I_{A\cap B})\}
\end{multline*} 
and
\begin{equation*}
	\mathscr{C}_2=\{C\in\mathscr{F}^{(1)}(s_1)\vee\mathscr{F}^{(2)}(s_2):\mathbb{E}(M_1(t_1)M_2(t_2)I_{C})=\mathbb{E}(M_1(s_1)M_2(s_2)I_{C})\}.
\end{equation*}	
Then, $\mathscr{C}_1$ and $\mathscr{C}_2$ are $\pi$-system and $\lambda$-system, respectively, such that $\mathscr{C}_1\subseteq\mathscr{C}_2$. Again, using Dynkin's $\pi$-$\lambda$ theorem, we have $\mathscr{C}_2=\mathscr{F}^{(1)}(t_1)\vee\mathscr{F}^{(2)}(t_2)$. Thus, $\{M_1(t_1)M_2(t_2),\ (t_1,t_2)\in\mathbb{R}^2_+\}$ is  $\mathscr{F}^{(1)}(t_1)\vee\mathscr{F}^{(2)}(t_2)$-martingale. This completes the proof.
\subsection{Proof of Theorem \ref{mthm3}}\label{app2}
By define $\mathcal{C}_1$, $\mathcal{C}_2$, $\mathscr{C}_1$ and $\mathscr{C}_2$ as follows:
\begin{align*}
	\mathcal{C}_1&=\{A\cap B\in\mathscr{F}_1(\textbf{s})\vee\mathscr{F}_2(\textbf{s}):\, A\in\mathscr{F}_1(\textbf{s}),\,B\in\mathscr{F}_2(\textbf{s})\\ &\hspace{2.5cm} \text{and}\ \mathbb{E}((M_1(\textbf{t})+M_2(\textbf{t}))I_{A\cap B})=\mathbb{E}((M_1(\textbf{s})+M_2(\textbf{s}))I_{A\cap B})\},\\
	\mathcal{C}_2&=\{C\in\mathscr{F}_1(\textbf{s})\vee\mathscr{F}_2(\textbf{s}):\mathbb{E}((M_1(\textbf{t})+M_2(\textbf{t}))I_{C})=\mathbb{E}((M_1(\textbf{s})+M_2(\textbf{s}))I_{C})\},\\
	\mathscr{C}_1&=\{A\cap B\in\mathscr{F}_1(\textbf{s})\vee\mathscr{F}_2(\textbf{s}):\, A\in\mathscr{F}_2(\textbf{s}),\,B\in\mathscr{F}_2(\textbf{s})\\ &\hspace{2.5cm} \text{and}\ \mathbb{E}(M_1(\textbf{t})M_2(\textbf{t})I_{A\cap B})=\mathbb{E}(M_1(\textbf{s})M_2(\textbf{s})I_{A\cap B})\}
\end{align*}	
and
\begin{equation*}
	\mathscr{C}_2=\{C\in\mathscr{F}_1(\textbf{s})\vee\mathscr{F}(\textbf{s}):\mathbb{E}(M_1(\textbf{t})M_2(\textbf{t})I_{C})=\mathbb{E}(M_1(\textbf{s})M_2(\textbf{s})I_{C})\},
\end{equation*}	
the proof follows similar lines to that of Theorem \ref{mthm1}. Hence, we omit it.
\subsection{Proof of Theorem \ref{mthm2}}\label{app3}
Suppose (i) holds, and $0\leq s_i\leq t_i$ for all $i=1,2,\dots,d$. For any $1\leq d'\leq d$, using Theorem 34.4 of \cite{Billingsley1995}, we get
\begin{align*}
	\mathbb{E}&\bigg(\sum_{i=1}^{d'}M_i(t_i)\Big|\vee_{i=1}^{d'}\mathscr{F}^{(i)}(s_i)\bigg)\\
	&=\mathbb{E}\bigg(\mathbb{E}\bigg(\sum_{i=1}^{d'}M_i(t_i)\Big|\vee_{j=1}^{d}\mathscr{F}^{(j)}(s_j)\bigg)\Big|\vee_{i=1}^{d'}\mathscr{F}^{(i)}(s_i)\bigg)\\
	&=\mathbb{E}\bigg(\sum_{j=1}^{d}M_j(s_j)-\mathbb{E}\bigg(\sum_{j=d'+1}^{d}M_j(t_j)\Big|\vee_{j=1}^{d}\mathscr{F}^{(j)}(s_j)\bigg)\Big|\vee_{i=1}^{d'}\mathscr{F}^{(i)}(s_i)\bigg)\\
	&=\mathbb{E}\bigg(\sum_{j=1}^{d}M_j(s_j)\Big|\vee_{i=1}^{d'}\mathscr{F}^{(i)}(s_i)\bigg)-\mathbb{E}\bigg(\sum_{j=d'+1}^{d}M_j(t_j)\Big|\vee_{i=1}^{d'}\mathscr{F}^{(i)}(s_i)\bigg)\\
	&=\sum_{i=1}^{d'}M_i(s_i)+\sum_{j=d'+1}^{d}\mathbb{E}(M_j(s_j)-M_j(t_j))=\sum_{i=1}^{d'}M_i(s_i),
\end{align*}
where we have used the independence of $\{M_j(t),\ t\ge0\}$'s to obtain the penultimate step. This proves the implication (i)$\implies$(ii).

For (ii)$\implies$(iii), it is sufficient to establish the martingale condition. To show this, we use induction. For $d=1$, the result holds true. Let us assume that the result holds for some $d=m\ge2$. If $\sum_{i=1}^{d'}M_i(t_i)$ is a $d'$-parameters $\vee_{i=1}^{d'}\mathscr{F}^{(i)}(t_i)$-martingale for $1\leq d'\leq m+1$ then for $s_i\leq t_i$, $i=1,2,\dots,m+1$, using Theorem 34.4 of \cite{Billingsley1995}, we have
\begin{align}
	\mathbb{E}\bigg(\sum_{i=1}^{m+1}M_i(t_i)\Big|\mathscr{F}^{(m+1)}(s_{m+1})\bigg)
	&=\mathbb{E}\bigg(\mathbb{E}\bigg(\sum_{i=1}^{m+1}M_i(t_i)\Big|\vee_{i=1}^{m+1}\mathscr{F}^{(i)}(s_i)\bigg)\Big|\mathscr{F}^{(m+1)}(s_{m+1})\bigg)\nonumber\\
	&=\mathbb{E}\bigg(\sum_{i=1}^{m+1}M_i(s_i)\Big|\mathscr{F}^{(m+1)}(s_{m+1})\bigg)\nonumber\\
	&=\sum_{i=1}^{m}\mathbb{E}M_i(s_i)+M_{m+1}(s_{m+1}),\label{smgp1}
\end{align}
where the last equality follows from the independence of $\{M_i(t_i),\ t_i\ge0\}$'s. Moreover,
\begin{align}
	\mathbb{E}\bigg(\sum_{i=1}^{m+1}M_i(t_i)\Big|\mathscr{F}^{(m+1)}(s_{m+1})\bigg)
	&=\sum_{i=1}^{m}\Big(M_i(t_i)|\mathscr{F}^{(m+1)}(s_{m+1})\Big)+\mathbb{E}(M_{m+1}(t_{m+1})|\mathscr{F}^{(m+1)}(s_{m+1}))\nonumber\\
	&=\sum_{i=1}^{m}\mathbb{E}M_i(t_i)+\mathbb{E}(M_{m+1}(t_{m+1})|\mathscr{F}^{(m+1)}(s_{m+1})).\label{smgp2}
\end{align}
From the induction hypothesis, it follows that $\mathbb{E}M_i(t_i)=\mathbb{E}M_i(s_i)$ for each $i=1,2,\dots,m$. On equating the right hand side of (\ref{smgp1}) to that of (\ref{smgp2}), we get 
\begin{equation*}
	\mathbb{E}(M_{m+1}(t_{m+1})|\mathscr{F}^{(m+1)}(s_{m+1}))=M_{m+1}(s_{m+1}),\ s_{m+1}\leq t_{m+1}.
\end{equation*} 
Thus, $\{M_{m+1}(t_{m+1}),\ t_{m+1}\ge0\}$ is a one parameter martingale with respect to its natural filtration.

(iii)$\implies$(i) follows from Theorem \ref{mthm1} and using the martingale property of $\{M_i(t_i),\ t_i\ge0\}$'s. This completes the proof of Theorem \ref{mthm2}.

\subsection{Proof of Lemma \ref{rcts}}\label{app4}
If $f_1,f_2,\dots,f_d$ are right continuous functions then for $\textbf{t}=(t_1,t_2,\dots,t_d)\in\mathbb{R}^d_+$, we get
\begin{equation*}
	\lim_{\substack{\textbf{s}\to\textbf{t}\\\textbf{t}\preceq\textbf{s}}}f(\textbf{s})=\lim_{\substack{\textbf{s}\to\textbf{t}\\\textbf{t}\preceq\textbf{s}}}\sum_{i=1}^{d}f_i(s_i)=\sum_{i=1}^{d}\lim_{s_i\to t_i^+}f_i(s_i)=\sum_{i=1}^{d}f_i(t_i)=f(\textbf{t}).
\end{equation*}
Thus, $f$ is a multiparameter right continuous function.

Let $a_1, a_2,\dots,a_d$ be the points of right continuity for $f_1$, $f_2$, $\dots$, $f_d$, respectively, and let $\textbf{a}(i)=(a_1,a_2,\dots,a_{i-1},t_i$, $a_{i+1},\dots,a_d)$ for all $i=1,2,\dots,d$. If $f$ is right continuous then  
\begin{equation*}
	\sum_{j\ne i}^{d}f_j(a_j)+f_i(t_i)=\lim_{\substack{\textbf{s}\to\textbf{a}(i)\\\textbf{a}(i)\preceq\textbf{s}}}\sum_{i=1}^{d}f_i(s_i)=\sum_{j\ne i}^{d}f_j(a_j)+\lim_{s_i\to t_i^+}f_i(s_i).
\end{equation*}
Therefore, $f_i$ is a right continuous function for each $i=1,2,\dots,d$. This completes the proof.

\end{document}